% With \loadreferencesfalse, the file makes its own counters and saves them to "references.tex".
% With \loadreferencestrue, the file loads the counters saved in "references.tex".
%
\newif\ifloadreferences\loadreferencestrue
%
%%%%%%%%%%%%%%%%%%%%%%%%%%%%%%%%%%%%%%%%%%%%%%%%%%%%%%%%%%%%%%%%%%%%%%%%%%%%%%%%%%%%%%%%%%%%%%%%%%%%%%%%%%%%%%%%%%%%%%%
%
% 0: Preliminaries.
%
%%%%%%%%%%%%%%%%%%%%%%%%%%%%%%%%%%%%%%%%%%%%%%%%%%%%%%%%%%%%%%%%%%%%%%%%%%%%%%%%%%%%%%%%%%%%%%%%%%%%%%%%%%%%%%%%%%%%%%%
%
%
\let\myfrac=\frac%
\input eplain %
\let\frac=\myfrac%
\let\myfootnote=\footnote%
\input amstex \input epsf %
\let\footnote=\myfootnote%
%
% Here we load the functions permitting us to use "amsmath" without using "amsppt".
%
\loadeufm\loadmsam\loadmsbm\message{symbol names}\UseAMSsymbols\message{,}%
\magnification 1200 %
\font\myfontdefault=cmr10%
\newif\ifmakebiblio%
\newif\ifinappendices%
\newif\ifundefinedreferences%
\newif\ifchangedreferences%
\makebibliofalse%
\undefinedreferencesfalse%
\changedreferencesfalse%
%
%%%%%%%%%%%%%%%%%%%%%%%%%%%%%%%%%%%%%%%%%%%%%%%%%%%%%%%%%%%%%%%%%%%%%%%%%%%%%%%%%%%%%%%%%%%%%%%%%%%%%%%%%%%%%%%%%%%%%%%
%
% 1: Abstract machinery.
%
% Here we define the macro "makecounter", which gives functionality to the counters defined presently. In order
% to implement it, it is necessary to provisionally change the catcodes.
%
%%%%%%%%%%%%%%%%%%%%%%%%%%%%%%%%%%%%%%%%%%%%%%%%%%%%%%%%%%%%%%%%%%%%%%%%%%%%%%%%%%%%%%%%%%%%%%%%%%%%%%%%%%%%%%%%%%%%%%%
%
\def\setcatcodes{\catcode`\!=0 \catcode`\\=11}%
{\global\let\noe=\noexpand%
\catcode`\@=11 \catcode`\_=11 \setcatcodes%
!global!def!_@@internal@@makeref#1{%
!global!expandafter!def!csname #1ref!endcsname##1{%
!csname _@#1@##1!endcsname%
!expandafter!ifx!csname _@#1@##1!endcsname!relax%
    !write16{#1 ##1 not defined - run saving references}%
    !undefinedreferencestrue%
!fi}}%
!global!def!_@@internal@@makelabel#1{%
!global!expandafter!def!csname #1label!endcsname##1{%
!edef!temptoken{!csname #1info!endcsname}%
!ifloadreferences%
!expandafter!ifx!csname _@#1@##1!endcsname!relax%
!write16{#1 ##1 not hitherto defined - rerun saving references}%
!changedreferencestrue%
!else%
!expandafter!ifx!csname _@#1@##1!endcsname!temptoken%
!else%
!write16{#1 ##1 reference has changed - rerun saving references}%
!changedreferencestrue%
!fi%
!fi%
!else%
!expandafter!edef!csname _@#1@##1!endcsname{!temptoken}%
!edef!textoutput{!write!references{\global\def\_@#1@##1{!temptoken}}}%
!textoutput%
!fi}}%
!global!def!makecounter#1{!_@@internal@@makelabel{#1}!_@@internal@@makeref{#1}}%
!unsetcatcodes%
}
%
%%%%%%%%%%%%%%%%%%%%%%%%%%%%%%%%%%%%%%%%%%%%%%%%%%%%%%%%%%%%%%%%%%%%%%%%%%%%%%%%%%%%%%%%%%%%%%%%%%%%%%%%%%%%%%%%%%%%%%%
%
% 2: Counters.
%
% Here we define the various counters.
%
%%%%%%%%%%%%%%%%%%%%%%%%%%%%%%%%%%%%%%%%%%%%%%%%%%%%%%%%%%%%%%%%%%%%%%%%%%%%%%%%%%%%%%%%%%%%%%%%%%%%%%%%%%%%%%%%%%%%%%%
%
\def\turnintolatin#1{\ifcase #1 _\or i\or ii\or iii\or iv\or v\or vi\or vii\or viii\or ix\or x\or xi\or xii\or xiii\or xiv\or xv\or xvi\or xvii\or xviii\or xix\or xx\or xxi\or xxii\or xxiii\or xxiv\or xxv\or xxvi\fi}%
\def\alphanum#1{\ifcase #1 _\or A\or B\or C\or D\or E\or F\or G\or H\or I\or J\or K\or L\or M\or N\or O\or P\or Q\or R\or S\or T\or U\or V\or W\or X\or Y\or Z\fi}%
\newwrite\references%
\ifloadreferences{\catcode`\@=11 \catcode`\_=11 \global\def\_@citation@Ahlfors{1}
\global\def\_@citation@AlcaldeCuestaDalBoMartinezVerjovsky{2}
\global\def\_@citation@AlessandriniLiuPapadopoulosSuI{3}
\global\def\_@citation@AlessandriniLiuPapadopoulosSuII{4}
\global\def\_@citation@AlessandriniLiuPapadopoulosSuIII{5}
\global\def\_@citation@AlvarezBrum{6}
\global\def\_@citation@AlvarezBrumMartinezPotrie{7}
\global\def\_@citation@AlvarezLessa{8}
\global\def\_@citation@BallmanSchroederGromov{9}
\global\def\_@citation@CandelI{10}
\global\def\_@citation@CandelConlon{11}
\global\def\_@citation@Chavel{12}
\global\def\_@citation@ChowKnopf{13}
\global\def\_@citation@Deroin{14}
\global\def\_@citation@DumasWolf{15}
\global\def\_@citation@EpsteinMillettTischler{16}
\global\def\_@citation@GhysI{17}
\global\def\_@citation@GilbTrud{18}
\global\def\_@citation@Hector{19}
\global\def\_@citation@HirschPughShub{20}
\global\def\_@citation@KatokHasselblatt{21}
\global\def\_@citation@LyubichMinsky{22}
\global\def\_@citation@McMullen{23}
\global\def\_@citation@MooreSchochet{24}
\global\def\_@citation@Rosenberg{25}
\global\def\_@citation@Saric{26}
\global\def\_@citation@ScannellWolf{27}
\global\def\_@citation@Spivak{28}
\global\def\_@citation@SullivanI{29}
\global\def\_@citation@SullivanII{30}
\global\def\_@citation@Tromba{31}
\global\def\_@citation@Verjovsky{32}
\global\def\_@citation@Wolpert{33}
\global\def\_@head@Introduction{1}
\global\def\_@subhead@TeichmuellerTheoryOfLaminations{1.1}
\global\def\_@subhead@PreviousResults{1.2}
\global\def\_@proc@Deroin{1.1}
\global\def\_@subhead@MainResults{1.3}
\global\def\_@proc@MainTheoremA{1.2}
\global\def\_@proc@MainTheoremB{1.3}
\global\def\_@subhead@GraftingsAndEarthquakesIntro{1.4}
\global\def\_@eqn@VariationOfLengthWithGraftings{\relax \unhbox \voidb@x \hbox {{\relax \tenrm (1)}}}
\global\def\_@proc@MainResultC{1.4}
\global\def\_@head@AnalyticPreliminaries{2}
\global\def\_@subhead@GreensFunctionsWithPointSingularities{2.1}
\global\def\_@eqn@LinearisedCurvatureOperator{\relax \unhbox \voidb@x \hbox {{\relax \tenrm (2)}}}
\global\def\_@eqn@CurvatureOperator{\relax \unhbox \voidb@x \hbox {{\relax \tenrm (3)}}}
\global\def\_@eqn@GreensFunctionCondition{\relax \unhbox \voidb@x \hbox {{\relax \tenrm (4)}}}
\global\def\_@eqn@GreensFunctionA{\relax \unhbox \voidb@x \hbox {{\relax \tenrm (5)}}}
\global\def\_@proc@GeneralGreensFunctionI{2.1}
\global\def\_@eqn@DecayRateOfPointGreensFunction{\relax \unhbox \voidb@x \hbox {{\relax \tenrm (6)}}}
\global\def\_@proc@PropertiesOfGreensOperatorInHyperbolicSpace{2.2}
\global\def\_@eqn@DecayAwayFromSource{\relax \unhbox \voidb@x \hbox {{\relax \tenrm (7)}}}
\global\def\_@proc@ExponentialDecayAwayFromSource{2.3}
\global\def\_@subhead@GreensFunctionsWithGeodesicSingularities{2.2}
\global\def\_@eqn@GreensFunctionConditionII{\relax \unhbox \voidb@x \hbox {{\relax \tenrm (8)}}}
\global\def\_@proc@FirstFormulaForGreensFunction{2.4}
\global\def\_@eqn@GreensFunctionB{\relax \unhbox \voidb@x \hbox {{\relax \tenrm (9)}}}
\global\def\_@proc@FormulaForGreensFunctionInClepsydra{2.5}
\global\def\_@eqn@DecayRateOfGeodesicGreensFunction{\relax \unhbox \voidb@x \hbox {{\relax \tenrm (10)}}}
\global\def\_@proc@PropertiesOfGreensFunctionGeodesic{2.6}
\global\def\_@eqn@GeneralGeodesicGreensFunction{\relax \unhbox \voidb@x \hbox {{\relax \tenrm (11)}}}
\global\def\_@proc@GeneralGreensFunctionGeodesic{2.7}
\global\def\_@proc@PropertiesOfGeodesicGreensFunction{2.8}
\global\def\_@eqn@DecayEstimateGeodesicGreensFunction{\relax \unhbox \voidb@x \hbox {{\relax \tenrm (12)}}}
\global\def\_@proc@PropertiesOfGeodesicGreensFunctionII{2.9}
\global\def\_@head@GraftingsAndEarthquakes{3}
\global\def\_@subhead@GraftingsAndEarthquakes{3.1}
\global\def\_@eqn@RatesOfChangeOfLengthI{\relax \unhbox \voidb@x \hbox {{\relax \tenrm (13)}}}
\global\def\_@eqn@RatesOfChangeOfLengthII{\relax \unhbox \voidb@x \hbox {{\relax \tenrm (14)}}}
\global\def\_@proc@RatesOfChangeOfLength{3.1}
\global\def\_@eqn@TwistedMetric{\relax \unhbox \voidb@x \hbox {{\relax \tenrm (15)}}}
\global\def\_@eqn@PerturbedMetric{\relax \unhbox \voidb@x \hbox {{\relax \tenrm (16)}}}
\global\def\_@eqn@SizeOfGrafting{\relax \unhbox \voidb@x \hbox {{\relax \tenrm (17)}}}
\global\def\_@eqn@SizeOfEarthquake{\relax \unhbox \voidb@x \hbox {{\relax \tenrm (18)}}}
\global\def\_@proc@ModulesOfSmoothTransformations{3.2}
\global\def\_@subhead@FirstOrderVariationsOfCurvatureAndLength{3.2}
\global\def\_@eqn@ThreeParameterPerturbedMetric{\relax \unhbox \voidb@x \hbox {{\relax \tenrm (19)}}}
\global\def\_@eqn@ConformalChangeOfCurvature{\relax \unhbox \voidb@x \hbox {{\relax \tenrm (20)}}}
\global\def\_@eqn@DerivativeOfCurvatureWithRespectToThirdComponent{\relax \unhbox \voidb@x \hbox {{\relax \tenrm (21)}}}
\global\def\_@eqn@ElementaryFormulaForDerivativeOfLength{\relax \unhbox \voidb@x \hbox {{\relax \tenrm (22)}}}
\global\def\_@proc@ElementaryFormulaForDerivativeOfLength{3.4}
\global\def\_@eqn@PerturbedCurvatureFormula{\relax \unhbox \voidb@x \hbox {{\relax \tenrm (23)}}}
\global\def\_@proc@PerturbedCurvatureFormula{3.5}
\global\def\_@eqn@ConnectionFormFormula{\relax \unhbox \voidb@x \hbox {{\relax \tenrm (24)}}}
\global\def\_@eqn@MovingFrameCurvature{\relax \unhbox \voidb@x \hbox {{\relax \tenrm (25)}}}
\global\def\_@eqn@DerivativeOfCurvature{\relax \unhbox \voidb@x \hbox {{\relax \tenrm (26)}}}
\global\def\_@eqn@DistributionalLimitOfDerivativeOfCurvature{\relax \unhbox \voidb@x \hbox {{\relax \tenrm (27)}}}
\global\def\_@proc@GaussBonnet{3.6}
\global\def\_@eqn@GraftingVariationFormula{\relax \unhbox \voidb@x \hbox {{\relax \tenrm (28)}}}
\global\def\_@eqn@EarthquakeVariationFormula{\relax \unhbox \voidb@x \hbox {{\relax \tenrm (29)}}}
\global\def\_@head@Laminations{4}
\global\def\_@subhead@Laminations{4.1}
\global\def\_@proc@PartitionOfUnity{4.1}
\global\def\_@eqn@LeafwiseBanachNorm{\relax \unhbox \voidb@x \hbox {{\relax \tenrm (30)}}}
\global\def\_@proc@Candel{4.2}
\global\def\_@proc@CandelSmoothness{4.3}
\global\def\_@subhead@TheDifferentialTopologyOfLaminations{4.2}
\global\def\_@eqn@DefinitionOfLeafwiseDistance{\relax \unhbox \voidb@x \hbox {{\relax \tenrm (31)}}}
\global\def\_@proc@SepMetComp{4.4}
\global\def\_@proc@BoundedAreaAndDiameter{4.5}
\global\def\_@proc@SinglePlaque{4.6}
\global\def\_@proc@TubularNeighbourhood{4.7}
\global\def\_@proc@LocalExtensionOfOpenSets{4.8}
\global\def\_@proc@ExtensionOfBumpFunctions{4.9}
\global\def\_@subhead@HyperbolicPerturbations{4.3}
\global\def\_@eqn@VariationOfMetricsFromA{\relax \unhbox \voidb@x \hbox {{\relax \tenrm (32)}}}
\global\def\_@eqn@HyperbolicPerturbationDefinition{\relax \unhbox \voidb@x \hbox {{\relax \tenrm (33)}}}
\global\def\_@proc@FormulaForDelta{4.10}
\global\def\_@eqn@FirstDefinitionOfDivergence{\relax \unhbox \voidb@x \hbox {{\relax \tenrm (34)}}}
\global\def\_@eqn@SecondDefintionOfDivergence{\relax \unhbox \voidb@x \hbox {{\relax \tenrm (35)}}}
\global\def\_@eqn@HyperbolicPerturbationCondition{\relax \unhbox \voidb@x \hbox {{\relax \tenrm (36)}}}
\global\def\_@proc@HyperbolicPerturbationCondition{4.11}
\global\def\_@subhead@InfiniteDimensionalityOfTeichmuellerSpace{4.4}
\global\def\_@proc@InfiniteDimensionality{4.12}
\global\def\_@proc@InductionStep{4.13}
\global\def\_@proc@IntegrabilityOfHyperbolicPerturbations{4.14}
\global\def\_@head@Bibliography{5}
 }%
\else{\openout\references=references.tex }%
\fi%
%
% A: Headings.
%
\newcount\headno%
\global\headno=0%
\def\headinfo{\ifinappendices\alphanum\headno\else\the\headno\fi}%
\def\nextheadno{\global\advance\headno by 1 \global\subheadno=0 \global\procno=0 \headinfo}%
\makecounter{head}%
%
% B: Subheadings.
%
\newcount\subheadno%
\global\subheadno=0%
\def\subheadinfo{\headinfo.\the\subheadno}%
\def\nextsubheadno{\global\advance\subheadno by 1 \subheadinfo}%
\makecounter{subhead}%
%
% C: Proclaims (Theorems, Propositions, Lemmas, Corollories, Definitions).
%
\newcount\procno%
\global\procno=0%
\def\procinfo{\headinfo.\the\procno}%
\def\nextprocno{\global\advance\procno by 1 \procinfo}%
\makecounter{proc}%
%
% D: Figures.
%
\newcount\figno%
\global\figno=0%
\def\figinfo{\subheadinfo.\the\figno}%
\def\nextfigno{\global\advance\figno by 1 \figinfo}%
\makecounter{fig}%
%
% E: Equations.
%
\newcount\eqnno%
\global\eqnno=0%
\def\eqninfo{\text{{\rm (\the\eqnno)}}}%
\def\nexteqnno[#1]{\global\advance\eqnno by 1 \eqninfo\hbox{\eqnlabel{#1}}}%
\makecounter{eqn}%
%
%%%%%%%%%%%%%%%%%%%%%%%%%%%%%%%%%%%%%%%%%%%%%%%%%%%%%%%%%%%%%%%%%%%%%%%%%%%%%%%%%%%%%%%%%%%%%%%%%%%%%%%%%%%%%%%%%%%%%%%
%
% 3: Citations.
%
% Citations are treated as a special type of counter.
%
%%%%%%%%%%%%%%%%%%%%%%%%%%%%%%%%%%%%%%%%%%%%%%%%%%%%%%%%%%%%%%%%%%%%%%%%%%%%%%%%%%%%%%%%%%%%%%%%%%%%%%%%%%%%%%%%%%%%%%%
%
\def\gobbleeight#1#2#3#4#5#6#7#8{}%
\newcount\citationno%
\global\citationno=0%
\def\citationinfo{\the\citationno}%
\makecounter{citation}%
\newwrite\biblio%
\def\newref#1#2{%
\def\temptext{#2}%
\edef\bibliotextoutput{\expandafter\gobbleeight\meaning\temptext}%
\global\advance\citationno by 1\citationlabel{#1}%
\ifmakebiblio%
    \edef\fileoutput{\write\biblio{\noindent\hbox to 0pt{\hss$[\the\citationno]$}\hskip 0.2em\bibliotextoutput\medskip}}%
    \fileoutput%
\fi}%
\def\cite#1{%
$[\citationref{#1}]$%
\ifmakebiblio%
    \edef\fileoutput{\write\biblio{#1}}%
    \fileoutput%
\fi%
}%
%
%%%%%%%%%%%%%%%%%%%%%%%%%%%%%%%%%%%%%%%%%%%%%%%%%%%%%%%%%%%%%%%%%%%%%%%%%%%%%%%%%%%%%%%%%%%%%%%%%%%%%%%%%%%%%%%%%%%%%%%
%
% 4: Formatting.
%
%%%%%%%%%%%%%%%%%%%%%%%%%%%%%%%%%%%%%%%%%%%%%%%%%%%%%%%%%%%%%%%%%%%%%%%%%%%%%%%%%%%%%%%%%%%%%%%%%%%%%%%%%%%%%%%%%%%%%%%
%
\let\mypar=\par%
\edef\Pagetitle={Blank}\headline={\hfil\Pagetitle\hfil}%
\edef\Pagefooter={Blank}\footline={\hfil\Pagefooter\hfil}%
%
% A: Move to next odd page.
%
\newcount\showpagenumflag%
\global\showpagenumflag=0 %
\def\nextoddpage%
{\newpage\ifodd\pageno%
\else\global\showpagenumflag=0 %
\null\vfil\eject%
\global\showpagenumflag=1 %
\fi}%
%
% B: Headings.
%
\font\headfont=cmb12%
\def\newhead#1[#2]%
{\ifhmode\mypar\fi%
\ifnum\headno=0 \else\goodbreak\bigskip\fi%
{\headfont\noindent\nextheadno\ - #1.}\headlabel{#2}%
\nobreak\medskip}%
%
% C: Subheadings.
%
\def\newsubhead#1[#2]%
{\ifhmode\mypar\fi%
\ifnum\subheadno=0 \else\goodbreak\medskip\fi%
{\bf\noindent\nextsubheadno\ - #1.\ }\subheadlabel{#2}}%
%
% D: Proclaims.
%
\newif\ifinproclaim%
\global\inproclaimfalse%
\def\proclaim#1{%
\goodbreak\medskip
\bgroup\inproclaimtrue%
\noindent{\bf #1}%
\nobreak\medskip\sl}%
\def\noskipproclaim#1{%
\goodbreak\medskip%
\bgroup\inproclaimtrue%
\noindent{\bf #1}\nobreak\sl}%
\def\endproclaim{\mypar\egroup\nobreak\medskip\ignorespaces}%
%
% E: Figures.
%
% The macro "\makelabelgrid" is a useful utility for guiding the positioning of labels is figures.
%
\newcount\xpos\newcount\ypos
\def\makelabelgrid{%
\xpos=-5 \ypos=-5 %
\loop\ifnum\xpos<6 %
{\loop\ifnum\ypos<6 %
\def\labeltext{x}%
\ifnum\xpos=0\def\labeltext{+}\fi%
\ifnum\ypos=0\def\labeltext{+}\fi%
\placelabel[\xpos][\ypos]{\labeltext}%
\advance\ypos by 1 %
\repeat}%
\advance\xpos by 1 %
\repeat}%
\def\placelabel[#1][#2]#3{{%
\setbox10=\hbox{\raise #2cm \hbox{\hskip #1cm #3}}%
\ht10=0pt \dp10=0pt \wd10=0pt \box10}}%
%
%
%
% F: Items.
%
\def\myitem#1{\noindent\hbox to .5cm{\hfill#1\hss}}%
%
% G: Right justification.
%
%
%
%%%%%%%%%%%%%%%%%%%%%%%%%%%%%%%%%%%%%%%%%%%%%%%%%%%%%%%%%%%%%%%%%%%%%%%%%%%%%%%%%%%%%%%%%%%%%%%%%%%%%%%%%%%%%%%%%%%%%%%
%
% 5: Special fonts.
%
%%%%%%%%%%%%%%%%%%%%%%%%%%%%%%%%%%%%%%%%%%%%%%%%%%%%%%%%%%%%%%%%%%%%%%%%%%%%%%%%%%%%%%%%%%%%%%%%%%%%%%%%%%%%%%%%%%%%%%%
%
% A: The "mathsf" font is not defined in Plain.
%
%
\font\sansseriften=cmss10%
\font\sansserifseven=cmss7%
\font\sansseriffive=cmss5%
\newfam\sansseriffam%
\textfont\sansseriffam=\sansseriften%
\scriptfont\sansseriffam=\sansserifseven%
\scriptscriptfont\sansseriffam=\sansseriffive%
\def\mathsf{\fam\sansseriffam}%
%
%
% B: The "mathbf" font is not defined in Plain.
%
\font\boldten=cmb10%
\font\boldseven=cmb7%
\font\boldfive=cmb5%
\newfam\mathboldfam%
\textfont\mathboldfam=\boldten%
\scriptfont\mathboldfam=\boldseven%
\scriptscriptfont\mathboldfam=\boldfive%
\def\mathbf{\fam\mathboldfam}%
%
%
% C: Here we define the macros "mathi" and "mathj". This is not really necessary, since the macros "imath" and
% "jmath" perform the same function. Still, it makes for an interesting exercise.
%
\font\mycmmiten=cmmi10%
\font\mycmmiseven=cmmi7%
\font\mycmmifive=cmmi5%
\newfam\mycmmifam%
\textfont\mycmmifam=\mycmmiten%
\scriptfont\mycmmifam=\mycmmiseven%
\scriptscriptfont\mycmmifam=\mycmmifive%
\def\hexa#1{\ifcase #1 0\or 1\or 2\or 3\or 4\or 5\or 6\or 7\or 8\or 9\or A\or B\or C\or D\or E\or F\fi}%
\mathchardef\mathi="7\hexa\mycmmifam7B%
\mathchardef\mathj="7\hexa\mycmmifam7C%
%
% D: Here we define a few Hebrew letters: "mybeth", "mygimmel" and "mydaleth".
%
\font\mymsbmten=msbm10 at 8pt%
\font\mymsbmseven=msbm7 at 5.6pt%6
\font\mymsbmfive=msbm5 at 4pt%
\newfam\mymsbmfam%
\textfont\mymsbmfam=\mymsbmten%
\scriptfont\mymsbmfam=\mymsbmseven%
\scriptscriptfont\mymsbmfam=\mymsbmfive%
\mathchardef\mybeth="7\hexa\mymsbmfam69%
\mathchardef\mygimmel="7\hexa\mymsbmfam6A%
\mathchardef\mydaleth="7\hexa\mymsbmfam6B%
%
%%%%%%%%%%%%%%%%%%%%%%%%%%%%%%%%%%%%%%%%%%%%%%%%%%%%%%%%%%%%%%%%%%%%%%%%%%%%%%%%%%%%%%%%%%%%%%%%%%%%%%%%%%%%%%%%%%%%%%%
%
% 6: Mathematics operators and symbols.
%
%%%%%%%%%%%%%%%%%%%%%%%%%%%%%%%%%%%%%%%%%%%%%%%%%%%%%%%%%%%%%%%%%%%%%%%%%%%%%%%%%%%%%%%%%%%%%%%%%%%%%%%%%%%%%%%%%%%%%%%
%
\def\proof{{\noindent\bf Proof:\ }}%
\def\remark{{\noindent\bf Remark:\ }}%
\def\qed{~$\square$}%
\def\makeop#1{\global\expandafter\def\csname op#1\endcsname{{\text{#1}}}}%
\def\makeopsmall#1{\global\expandafter\def\csname op#1\endcsname{{\text{\lowercase{#1}}}}}%
%
% A: Set Theory.
%
\def\munion{\mathop{\cup}}%
\def\minter{\mathop{\cap}}%
%
% B: Point set topology.
%
\makeop{Ext}%
\makeop{Int}%
\makeop{Dist}%
\makeop{Diam}%
\makeop{Length}%
%
% C: Sequences.
%
%
%
%
\def\mlim{\mathop{{\text{Lim}}}}%
\def\msup{\mathop{{\text{Sup}}}}%
\def\minf{\mathop{{\text{Inf}}}}%
%
% D: Linear Algebra.
%
\makeop{Dim}%
\makeop{Ker}%
\makeop{Coker}%
\makeop{Tr}%
\makeop{Adj}%
\makeop{Det}%
\makeop{End}%
\makeop{Lin}%
\makeop{Symm}%
\makeop{Mult}%
%
% E: Basic calculus.
%
\makeop{dx}%
\makeop{dy}%
\makeop{dz}%
\makeop{dt}%
\makeop{dVol}%
\makeop{dArea}%
\makeop{Supp}%
\makeop{Hess}%
\makeop{Lip}%
%
% F: Complex Numbers.
%
\makeop{Re}%
\makeop{Im}%
\makeop{Arg}%
\makeop{Log}%
\makeop{Exp}%
%
% G: Trigonometry.
%
\makeopsmall{Cos}%
\makeopsmall{Sin}%
\makeopsmall{Tan}%
\makeopsmall{Sec}%
\makeopsmall{Cosec}%
\makeopsmall{Cot}%
\makeopsmall{ArcCos}%
\makeopsmall{ArcSin}%
\makeopsmall{ArcTan}%
\makeopsmall{ArcSec}%
\makeopsmall{ArcCosec}%
\makeopsmall{ArcCot}%
%
% H: Hyperbolic Trigonometry.
%
\makeopsmall{Cosh}%
\makeopsmall{Sinh}%
\makeopsmall{Tanh}%
\makeopsmall{ArcCosh}%
\makeopsmall{ArcSinh}%
\makeopsmall{ArcTanh}%
%
% I: Differential and Riemannian Geometry.
%
\makeop{Vol}%
\makeop{Area}%
\makeop{Riem}%
\makeop{Ric}%
\makeop{Scal}%
\makeop{Euc}%
\makeop{Imm}%
\makeop{Emb}%
%
% J: Lie Groups.
%
\makeop{Id}%
\makeop{Ad}%
\makeop{O}%
\makeop{SO}%
\makeop{SL}%
\makeop{GL}%
\makeop{Conf}%
\makeop{Homeo}%
\makeop{Diff}%
\makeop{Isom}%
%
% K: Functional Analysis.
%
\makeop{Ind}%
\makeop{Sig}%
\makeop{Spec}%
%
% L: Commutative Diagrams.
%
\def\harr#1#2{\smash{\mathop{\hbox to .5in{\rightarrowfill}}\limits^{\scriptstyle #1}_{\scriptstyle #2}}}%
%
%
%
% M: Other.
%
\makeop{Conv}%
\makeop{Max}%
\makeop{Min}%
\makeop{Mod}%
\makeop{Deg}%
\makeop{loc}%
%
%%%%%%%%%%%%%%%%%%%%%%%%%%%%%%%%%%%%%%%%%%%%%%%%%%%%%%%%%%%%%%%%%%%%%%%%%%%%%%%%%%%%%%%%%%%%%%%%%%%%%%%%%%%%%%%%%%%%%%%
%
% 7: Redundant Material.
%
%%%%%%%%%%%%%%%%%%%%%%%%%%%%%%%%%%%%%%%%%%%%%%%%%%%%%%%%%%%%%%%%%%%%%%%%%%%%%%%%%%%%%%%%%%%%%%%%%%%%%%%%%%%%%%%%%%%%%%%
%
% This file contains redundant functionality for constructing a bibliography. Although it is not used, it may one day
% prove useful, and so I leave it here.
%
% A: Before the citations, the file should contain:
%
% \newif\ifmakebiblio
%
% followed by either \makebibliotrue or \makebibliofalse.
%
% (note that, for conveniance, the commands \newif\ifmakebiblio and \makebibliofalse have been included at the
% beginning of this preamble. These commands should be removed before making the changes outlined here.
%
% B: Immediately before the first citation, the file should contain the following instructions:
%
% \ifmakebiblio%
% \openout\biblio=biblio.tex %
% \edef\fileoutput{\write\biblio{\bgroup\leftskip=2em}}%
% \fileoutput%
% \fi%
%
% C: Immediately after the last citation, the file should contain the following instruction:
%
% \ifmakebiblio%
% {\edef\fileoutput{\write\biblio{\egroup}}%
% \fileoutput}%
% \fi%
%
 %
%
% In order to use "dvips", enter:
% dvips -N0 -Z0 -K0 [whatever.dvi] -o [whatever.ps]
%
%%%%%%%%%%%%%%%%%%%%%%%%%%%%%%%%%%%%%%%%%%%%%%%%%%%%%%%%%%%%%%%%%%%%%%%%%%%%%%%%%%%%%%%%%%%%%%%%%%%%%%%%%%%%%%%%%%%%%%%
%
% 1: The Paper.
%
%%%%%%%%%%%%%%%%%%%%%%%%%%%%%%%%%%%%%%%%%%%%%%%%%%%%%%%%%%%%%%%%%%%%%%%%%%%%%%%%%%%%%%%%%%%%%%%%%%%%%%%%%%%%%%%%%%%%%%%
%
\font\tablefont=cmr7
\newif\ifshowaddress\showaddresstrue
\def\Pagetitle{\hfil}
\def\Pagefooter{\hfil}
\null \vfill
\def\centre{\rightskip=0pt plus 1fil \leftskip=0pt plus 1fil \spaceskip=.3333em \xspaceskip=.5em \parfillskip=0em \parindent=0em}%
\def\textmonth#1{\ifcase#1\or January\or Febuary\or March\or April\or May\or June\or July\or August\or September\or October\or November\or December\fi}
\font\abstracttitlefont=cmr10 at 14pt {\abstracttitlefont\centre Earthquakes and graftings of hyperbolic surface laminations.\par}
\bigskip
{\centre 26th July 2019\par}
%{\centre \the\day\ \textmonth\month\ \the\year\par}
\bigskip
{\centre S\'ebastien Alvarez\footnote{${}^1$}{{\tablefont CMAT, Facultad de Ciencias, Universidad de la República
Igua 4225 esq. Mataojo. Montevideo, Uruguay.\hfill}},
Graham Smith\footnote{${}^2$}{{\tablefont Instituto de Matem\'atica, UFRJ, Av. Athos da Silveira Ramos 149, Centro de Tecnologia - Bloco C, Cidade Universit\'aria - Ilha do Fund\~ao, Caixa Postal 68530, 21941-909, Rio de Janeiro, RJ - BRAZIL\hfill}}\par}
\bigskip
\noindent{\bf Abstract:~}We study compact hyperbolic surface laminations. These are a generalization of closed hyperbolic surfaces which appear to be more suited to the study of Teichm\"uller theory than arbitrary non-compact surfaces. We show that the Teichm\"uller space of any non-trivial hyperbolic surface lamination is infinite dimensional. In order to prove this result, we study the theory of deformations of hyperbolic surfaces, and we derive what we believe to be a new formula for the derivative of the length of a simple closed geodesic with respect to the action of grafting. This formula complements those derived by McMullen in \cite{McMullen}, in terms of the Weil-Petersson metric, and by Wolpert in \cite{Wolpert}, for the case of earthquakes.
\bigskip
\noindent{\bf Classification AMS~:~}30F60
%
% 30F60, Teichmueller Theory
% 53C50, Lorentz manifold
%
\par
\vfill
\eject
%\nextoddpage
%
\global\pageno=1
\myfontdefault
\def\Pagetitle{\hfil Earthquakes and graftings of hyperbolic surface lamintions.\hfil}
\def\Pagefooter{\hfil{\myfontdefault\folio}\hfil}
\catcode`\@=11
\def\triplealign#1{\null\,\vcenter{\openup1\jot \m@th %
\ialign{\strut\hfil$\displaystyle{##}\quad$&\hfil$\displaystyle{{}##}$&$\displaystyle{{}##}$\hfil\crcr#1\crcr}}\,}
\def\multiline#1{\null\,\vcenter{\openup1\jot \m@th %
\ialign{\strut$\displaystyle{##}$\hfil&$\displaystyle{{}##}$\hfil\crcr#1\crcr}}\,}
\catcode`\@=12
\newref{Ahlfors}{Ahlfors L. V., {\sl Conformal Invariants: Topics in Geometric Function Theory}, McGraw-Hill, New York, D\"usseldorf, Johannesburg, (1973)}
\newref{AlcaldeCuestaDalBoMartinezVerjovsky}{Alcalde Cuesta F., Dal'Bo F., Mart\'\i nez M., Verjovsky, Minimality of the horocycle flow
on laminations by hyperbolic surfaces with non-trivial topology, {\sl Discrete Contin. Dyn. Syst.}, {\bf 36}, no. 9, (2016), 4619--4635}
\newref{AlessandriniLiuPapadopoulosSuI}{Alessandrini A., Liu L., Papadopoulos A., Su W., Sun Z., On Fenchel-Nielsen coordinates
on Teichm\"ller spaces of surfaces of infinite type, {\sl Ann. Acad. Sci. Fenn. Math.}, {\bf 36}, no. 2, (2011), 621--659}
\newref{AlessandriniLiuPapadopoulosSuII}{Alessandrini A., Liu L., Papadopoulos A., Su W., On local comparison between various
metrics on Teichm\"uller spaces, {\sl Geom. Dedicata}, {\bf 157}, (2012), 91--110}
\newref{AlessandriniLiuPapadopoulosSuIII}{Alessandrini A., Liu L., Papadopoulos A., Su W., On various Teichm\"uller spaces of a
surface of infinite topological type, {\sl Proc. Amer. Math. Soc.}, {\bf 140}, no. 2, (2012), 561--574}
\newref{AlvarezBrum}{Alvarez S., Brum J., Topology of leaves for minimal laminations by hyperbolic surfaces II, Preprint, 2019}
\newref{AlvarezBrumMartinezPotrie}{Alvarez S., Brum J., Mart\'\i nez, Potrie R., Topology of leaves for minimal laminations by hyperbolic surfaces, with an appendix written with M. Wolff. Preprint, 2019}
\newref{AlvarezLessa}{Alvarez S., Lessa P., The Teichmüller space of the Hirsch foliation, {\sl Ann. Inst. Fourier}, {\bf 68}, no. 1, (2018), 1--51}
\newref{BallmanSchroederGromov}{Ballmann W., Gromov M., Schroeder V., {\sl Manifolds of Nonpositive Curvature}, Progress in Mathematics, {\bf 61}, Birkha\"user Verlag, (1985)}
\newref{CandelI}{Candel A., Uniformization of surface laminations, {\sl Ann. Sci. ENS.}, {\bf 4}, no. 26, (1993), 489--516}
\newref{CandelConlon}{Candel A., Conlon L., {\sl Foliations I}, Graduate Studies in Mathematics, {\bf 23}, AMS, Providence, Rhode Island, (2000)}
\newref{Chavel}{Chavel I., {\sl Riemannian Geometry: A Modern Introduction}, Cambridge Studies in Advanced Mathematics, {\bf 98}, CUP, (2006)}
\newref{ChowKnopf}{Chow B., Knopf D., {\sl The Ricci Flow: An Introduction}, Mathematical Surveys and Monographs, {\bf 110}, AMS, (2004)}
\newref{Deroin}{Deroin B., Nonrigidity of hyperbolic surfaces laminations, {\sl Proc. Amer. Math. Soc.}, {\bf 135}, no. 3, (2007), 873--881}
\newref{DumasWolf}{Dumas D., Wolf M., Projective structures, grafting and measured laminations, {\sl Geom. Topol.}, {\bf 12}, no. 1, (2008), 351--386}
\newref{EpsteinMillettTischler}{Epstein D. B. A., Millett K. C., Tischler D., Leaves without holonomy, {\sl J. London Math. Soc.}, {\bf 16}, no. 3, (1977), 548--552}
\newref{GhysI}{Ghys \'E., Laminations par surfaces de Riemann, in {\sl Dynamique et g\'eom\'trie complexes (Lyon,
1997)}, {\sl Panor. Synth\`eses}, Volume 8, Pages ix, xi, 49--95, Soc. Math. France, Paris, 1999}
\newref{GilbTrud}{Gilbarg D., Trudinger N., {\sl Elliptic partial differential equations of second order}, Grundlehren der Mathematischen Wissenschaften, {\bf 224}, Springer-Verlag, Berlin, second edition, 1983}
\newref{Hector}{Hector G., Feuilletages en cylindres, in {\sl Geometry and topology (Proc. III Latin Amer. School
of Math., Inst. Mat. Pura Aplicada CNPq, Rio de Janeiro, 1976)}, Volume 597 of {\sl Lecture Notes in Math.}, Pages 252--270, Springer, Berlin, (1977)}
\newref{HirschPughShub}{Hirsch M. W., Pugh C. C., Shub M., {\sl Invariant manifolds}, Lecture Notes in Mathematics, {\bf 583}, Springer Verlag, Berling, New York, (1977)}
\newref{KatokHasselblatt}{Katok A., Hasselblatt B., {\sl Introduction to the modern theory of dynamical systems}, Encylopedia of Mathematics and its Applications, {\bf 54}, CUP, Cambridge, (1995)}
\newref{LyubichMinsky}{Lyubich M., Minsky Y., Laminations in holomorphic dynamics, {\sl J. Diff. Geom.}, {\bf 47}, (1997), 17--94}
\newref{McMullen}{McMullen C., Complex earthquakes and Teichmüller theory, {\sl J. Amer. Math. Soc.}, {\bf 11}, no. 2, (1998), 283--320}
\newref{MooreSchochet}{Moore C. C., Schochet C. L., {\sl Global analysis on foliated spaces}, MSRI Publications, {\bf 9}, Cambridge University Press, New York, second edition, 2006}
\newref{Rosenberg}{Rosenberg H., Foliations by planes, {\sl Topology}, {\bf 7}, (1968), 131--138}
\newref{Saric}{\v Sari\'c, The Teichm\"uller theory of the solenoid, in {\sl Handbook of Teichmüller theory. Vol. II},
Volume 13 of {\sl IRMA Lect. Math. Theor. Phys.}, Pages 811--857, Eur. Math. Soc., Zürich, 2009}
\newref{ScannellWolf}{Scannell K., Wolf M., The grafting map of Teichmüller space, {\sl J. Amer. Math. Soc.}, {\bf 15}, no. 4, (2002), 893--927}
\newref{Spivak}{Spivak M., {\sl A Comprehensive Introduction to Differential Geometry}, Vols. 1-5, Publish of Perish, (1999)}\newref{SullivanI}{Sullivan D., Bounds, quadratic differentials, and renormalization conjectures, in {\sl American
Mathematical Society centennial publications, Vol. II, (Providence, RI, 1988)}, Pages 417--466, Amer. Math. Soc., Providence, RI, 1992}
\newref{SullivanII}{Sullivan D, Linking the universalities of Milnor-Thurston, Feigenbaum and Ahlfors-Bers, in {\sl Topological methods in modern mathematics (Stony Brook, NY, 1991)}, Pages 543--564, Publish or Perish, Houston, TX, 1993}
\newref{Tromba}{Tromba A., {\sl Teichm\"uller theory in Riemannian Geometry}, Lectures in Mathematics, ETH Z\"urich, Birkha\"user Verlag, 1992}
\newref{Verjovsky}{Verjovsky A., A uniformization theorem for holomorphic foliations, in {\sl The Lefschetz centennial conference, Part III (Mexico City, 1984)}, Volume 58 of {\sl Contemp. Math.}, Pages 233--253, Amer. Math. Soc., Providence, RI, 1987}
\newref{Wolpert}{Wolpert S., An elementary formula for the Fenchel-Nielsen twist, {\sl Comment. Math. Helv.}, {\bf 56}, no. 1, (1981), 132--135}
\makeopsmall{Sinh}%
\makeop{inj}%
\makeopsmall{Ln}%
\makeop{Conv}%
\makeop{GH}%
\makeopsmall{Coth}%
\makeop{hyp}%
\makeop{euc}%
\makeop{PSL}%
\makeop{M}%
\makeop{bdd}%
\makeop{dArea}%
\makeop{dl}%
\makeopsmall{Arccoth}%
\makeopsmall{Arccot}%
%
%%%%%%%%%%%%%%%%%%%%%%%%%%%%%%%%%%%%%%%%%%%%%%%%%%%%%%%%%%%%%%%%%%%%%%%%%%%%%%%%%%%%%%%%%%%%%%%%%%%%
%
%%%%%%%%%%%%%%%%%%%%%%%%%%%%%%%%%%%%%%%%%%%%%%%%%%%%%%%%%%%%%%%%%%%%%%%%%%%%%%%%%%%%%%%%%%%%%%%%%%%%
%
\makeop{Ann}%
\newhead{Introduction}[Introduction]
\newsubhead{Teichm\"uller theory of laminations}[TeichmuellerTheoryOfLaminations]
Laminations, which have various applications in the study of hyperbolic dynamics (c.f. \cite{HirschPughShub}, \cite{KatokHasselblatt}, \cite{LyubichMinsky} and \cite{SullivanI}), are an extension of foliations where the ambient space is no longer assumed to be a smooth manifold.\footnote*{A brief review of the fundamentals of the theory of laminations is given in Section \subheadref{Laminations}.} A hyperbolic surface lamination is a lamination in which all leaves are Riemann surfaces of hyperbolic type, that is, Riemann surfaces which are uniformized by the Poincar\'e disk. Such laminations arise quite frequently (see, for example, \cite{GhysI}), and the property of a given compact surface lamination being hyperbolic is {\sl topological} in the sense that it does not depend on the laminated conformal structure chosen (see \cite{CandelI}).
\par
Hyperbolic surface laminations also appear to possess a better structured Teichm\"uller theory than arbitrary non-compact hyperbolic surfaces. Indeed, in \cite{AlessandriniLiuPapadopoulosSuI}, \cite{AlessandriniLiuPapadopoulosSuII} and \cite{AlessandriniLiuPapadopoulosSuIII}, natural constructions of the Teichm\"uller space of a hyperbolic surface of infinite topological type are defined using pants decompositions, complex structures and length spectra. However, the authors then show that each of these different methods may yield a different space. On the other hand, the known natural constructions of the Teichm\"uller space of a compact hyperbolic surface lamination are all equivalent.
\par
In \cite{SullivanI} and \cite{SullivanII}, Sullivan defines the Teichm\"uller space of a hyperbolic surface lamination to be the set of transversally continuous conformal structures modulo leafwise diffeorphisms which are leafwise isotopic to the identity. An alternative description using leafwise complex structures is described by Moore \& Schochet in \cite{MooreSchochet}. In \cite{CandelI} (c.f. also \cite{Verjovsky}), Candel proves that every hyperbolic surface lamination carries a unique leafwise metric in its conformal class (see Section \subheadref{Laminations}) which allows the Teichm\"uller space of compact hyperbolic surface laminations to be studied as the space of hyperbolic leafwise metrics modulo leafwise diffeomorphisms which are leafwise isotopic to the identity. It is this framework that we will adopt in the sequel.
\newsubhead{Previous results}[PreviousResults]
Little is currently known about the general theory of the Teichm\"uller space of a given compact hyperbolic surface lamination. In \cite{SullivanII}, Sullivan shows that this space is a {\sl Banach manifold} carrying a natural complex structure with respect to which it is biholomorphic to an open subset of the space of leafwise holomorphic quadratic differentials. In \cite{Deroin}, motivated by Ghys' construction of non-constant meromorphic functions on hyperbolic laminations using Poincar\'e series (see \cite{GhysI}), Deroin proves
\proclaim{Theorem \nextprocno, Deroin \cite{Deroin}}
\noindent If a hyperbolic surface lamination contains a simply connected leaf, then its Teichm\"uller space is infinite dimensional.
\endproclaim
\proclabel{Deroin}
Beyond these general results, the authors are only aware of three specific cases in which the Teichm\"uller space of a compact hyperbolic surface lamination is understood. The first is that of the family of laminations, discussed in \cite{GhysI} and \cite{SullivanII}, associated to expanding maps of the unit circle. The second is that of Sullivan's universal solenoid, obtained as the inverse limit of finite coverings of a closed hyperbolic surface, whose Teichm\"uller space was computed by \v Sari\'c in \cite{Saric}. The third is that of the Hirsch foliation, defined as the quotient of the stable foliation of the hyperbolic attractor of Smale's solenoidal map, whose Teichm\"uller space was computed by the first author in collaboration with Lessa in \cite{AlvarezLessa}.
\newsubhead{Main results}[MainResults]
Our main result completes that of Deroin. A hyperbolic surface lamination will be said to be {\sl trivial} whenever it consists of a finite union of closed surfaces. We show
\proclaim{Theorem \nextprocno}
\noindent The Teichm\"uller space of a non-trivial hyperbolic surface lamination is infinite dimensional.
\endproclaim
\proclabel{MainTheoremA}
We will explain presently how Theorem \procref{MainTheoremA} follows immediately from Deroin's result and the second main result of this paper. Our approach takes advantage of the topology of the leaves to generate an arbitrarily large number of independent movements in Teichm\"uller space via perturbations of Candel's hyperbolic leafwise metric. In order to state the result, we introduce some notation. First, given a complete hyperbolic surface $\Sigma$ and a simple closed curve $\gamma$, let $[\gamma]$ denote the free homotopy class in which $\gamma$ lies and let $l([\gamma],g)$ denote the infimal length of curves in this class with respect to the metric $g$. Next, given a surface lamination $X$, its leafwise topology is defined to be the topology generated by all open subsets of leaves of this lamination (see Section \subheadref{TheDifferentialTopologyOfLaminations}). In particular, a subset $Y$ of $X$ is compact in this topology if and only if it is a finite union of compact subsets of leaves. We prove
\proclaim{Theorem \nextprocno}
\noindent Let $X$ be a compact hyperbolic surface lamination. Suppose that for every subset $Y$ of $X$ which is compact in the leafwise topology there exists a simple, closed leafwise geodesic $\gamma$ in $X$ not intersecting $Y$. Then there exists an infinite sequence $(\gamma_m)_{m\in\Bbb{N}}$ of simple, closed leafwise geodesics in $X$ such that, for every $m\in\Bbb{N}$ and for every finite sequence $a_1,...,a_m\in\Bbb{R}$, there exists a smooth family of leafwise hyperbolic metrics $(g_t)_{t\in]-\epsilon,\epsilon[}$ such that, for all $1\leq i\leq m$,
$$
\left.\frac{\partial}{\partial t}\opLog(l([\gamma_i],g_t))\right|_{t=0} = a_i.
$$
In particular, the Teichm\"uller space of $X$ is infinite dimensional.
\endproclaim
\proclabel{MainTheoremB}
\noindent Theorem \procref{MainTheoremA} follows from Theorems \procref{Deroin} and \procref{MainTheoremB} by an argument used in \cite{AlcaldeCuestaDalBoMartinezVerjovsky} and \cite{AlvarezBrum}. Indeed, it was proven independently by Epstein, Millett \& Tischler in \cite{EpsteinMillettTischler} and by Hector in \cite{Hector} that the generic leaf of a compact lamination has trivial holonomy. Consequently, if a hyperbolic surface lamination has no simply connected leaf, then there exists a simple, closed geodesic inside a leaf with trivial holonomy. By Reeb's stability theorem, this geodesic has a neighbourhood trivially laminated by annuli. It then follows by the transverse continuity of the leafwise metric and the persistence of closed geodesics under perturbations of hyperbolic metrics that each of these annuli also contains a simple, closed geodesic. Since the lamination is non-trivial, this yields sufficient simple, closed leafwise geodesics for the hypotheses of Theorem \procref{MainTheoremB} to be satisfied and Theorem \procref{MainTheoremA} follows.
\par
However, there are laminations that can be treated simultaneously by the methods of both theorems. Indeed, the laminations to which Deroin's theorem applies but not ours are precisely those for which all leaves are simply connected except for finitely many leaves which are of finite topological type. Although examples of such laminations were constructed in \cite{AlvarezBrumMartinezPotrie} for $3$-dimensional ambient spaces, we believe this condition to be quite restrictive, as this is the case for foliations. Indeed, in \cite{AlvarezBrumMartinezPotrie} it is shown that if all but a finite number of leaves a smooth, minimal foliation are simply connected, and if all the remaining leaves are of finite topological type, then all the leaves of the foliation are in fact planes. It then follows by the result \cite{Rosenberg} of Rosenberg that the ambient manifold is a $3$-dimensional torus and the foliation is by parabolic planes.
\newsubhead{Graftings and Earthquakes}[GraftingsAndEarthquakesIntro]
Our proof closely follows the ideas developed by the first author in collaboration with Lessa in \cite{AlvarezLessa}, where Fenchel-Nielsen coordinates were used to parametrize the Teichm\"uller space of the Hirsch foliation. In the general case, the existence of such coordinates cannot be guaranteed, and we thus make use of the curves given in the hypotheses of Theorem \procref{MainTheoremB} in order to define surgeries of the leafwise metric which vary {\sl independently} the lengths of an arbitrarily large number of curves.
\par
The surgeries that we use are generalisations of graftings, which we recall are defined as follows (c.f. \cite{DumasWolf}, \cite{McMullen} and \cite{ScannellWolf}). Given a simple, closed geodesic $\gamma$ in a marked, hyperbolic surface $\Sigma$ and a positive real number $t$, the {\sl grafting of length $t$} of $\Sigma$ along $\gamma$ is defined to be the marked, hyperbolic surface obtained by cutting $\Sigma$ along $\gamma$, inserting a flat cylinder of length $t$, and multiplying the metric of the resulting surface by a suitable conformal factor. A related surgery operation is that of earthquakes. For any real number $t$, the right earthquake of $\Sigma$ along $\gamma$ is defined to be the marked, hyperbolic surface obtained by rotating the right hand side of $\gamma$ a distance $t$ in the positive direction of this geodesic. The resulting marked, hyperbolic surfaces will be denote by $G(\gamma,\Sigma)(t)$ and $E(\gamma,\Sigma)(t)$ respectively.
\par
Our main result essentially consists in extending the grafting surgery to the framework of laminations, which we believe to be of independent interest. It seems unlikely that there exist a canonical way of extending graftings to laminations in general. Furthermore, although it is straightforward to extend to the entire lamination a grafting along a simple, closed geodesic with trivial holonomy, it is less clear how this can be done for geodesics with more general holonomy. Our construction yields extensions of graftings along arbitrary simple, closed geodesics. Furthermore, although these extensions are non-canonical, they can always be chosen so that their supports are contained in arbitrarily small neighbourhoods of the geodesic in question (c.f. Lemma \procref{InductionStep}).
\par
Theorem \procref{MainTheoremB} then follows by carefully estimating the effect of grafting on the lengths of all other simple, closed leafwise geodesics of the lamination. It is unsurprising and perfectly consistent with known results of hyperbolic geometry that the effect of a grafting should decay exponentially with distance from the locus of surgery. However, our calculations are simplified with the help of the following exact formula, which complements those derived by McMullen in \cite{McMullen} and Wolpert in \cite{Wolpert} (see also \cite{DumasWolf} and \cite{ScannellWolf}) and which, to our knowledge, has not previously appeared in the litterature. Using the notation introduced before the statement of Theorem \procref{MainTheoremB}, we show
\proclaim{Theorem \nextprocno}
\noindent For every pair $(\gamma,\gamma')$ of simple, closed geodesics in $\Sigma$,
\goodbreak
$$
\left.\frac{\partial}{\partial t}l([\gamma'],G(\gamma,\Sigma)(t))\right|_{t=0}
=\sum_{x\in\gamma\minter\gamma'}\opSin(\theta_x) + \int_\gamma\int_{\gamma'}K(p,q)dl_pdl_q.\eqnum{\nexteqnno[VariationOfLengthWithGraftings]}
$$
where $K(x,y)$ is the Green's kernel over $\Sigma$ of the operator
$$
L:=\Delta-2,
$$
and, for all $x\in\gamma\minter\gamma'$, $\theta_x$ denotes the angle that $\gamma$ makes with $\gamma'$ at the point $x$.
\endproclaim
\proclabel{MainResultC}
\noindent This result is proven in Theorem \procref{RatesOfChangeOfLength}, below.
\newsubhead{Acknowledgements}[Acknowledgements]
The authors are grateful to Thierry Barbot, David Dumas, Fran\c{c}ois Fillastre, Gabriel Calsamiglia and Pablo Lessa for helpful comments insightful conversations. The first author was partially supported by ANII via the Fondo Clemente Estable (projects FCE\_135352 and
FCE\_3\_2018\_1\_148740), by CSIC I+D 389, and by the MathAmSud project RGSD (Rigidity and Geometric Structures in Dynamics) 19-MATH-04. The second author was partially supported by the MathAmSud project GDAR (Geometry and Dynamics of Anosov Representations).
\newhead{Analytic preliminaries}[AnalyticPreliminaries]
\newsubhead{Green's functions with point singularities}[GreensFunctionsWithPointSingularities]
Let $\Sigma$ be a complete hyperbolic surface without cusps. Let $\Delta$ be its Laplace-Beltrami operator. Define
$$
L := \Delta - 2.\eqnum{\nexteqnno[LinearisedCurvatureOperator]}
$$
This operator is, up to a change of sign, the linearisation about the hyperbolic metric of the curvature operator (c.f. \cite{Chavel} and \cite{ChowKnopf})
$$
\kappa_\phi := -e^{-2\phi}\big(\Delta\phi - \kappa_0\big).\eqnum{\nexteqnno[CurvatureOperator]}
$$
It is common to study its properties using general elliptic theory. In the present context, however, it is simpler to determine explicit formulae for its inverses in terms of Green's functions. Recall first that, considered as an operator acting on distributions, $L$ has trivial kernel in the space $C^2_\opbdd(\Sigma)$ of bounded, twice-differentiable functions over $\Sigma$ as well as in the space $L^1(\Sigma)$ of integrable functions over $\Sigma$. For all $x\in\Sigma$, the {\sl Green's function} of $L$ over $\Sigma$ with singularity at $x$ is defined to be the unique function $K_x\in L^1(\Sigma)$ such that
$$
LK_x\opdArea = \delta_x,\eqnum{\nexteqnno[GreensFunctionCondition]}
$$
in the distributional sense, where $\opdArea$ denotes the area form of $\Sigma$, and $\delta_x$ denotes the {\sl Dirac delta distribution} with singularity at $x$, that is
$$
\langle\delta_x,g\rangle := f(x).
$$
The {\sl Green's kernel} of $L$ over $\Sigma$ is defined such that, for all $x\neq y\in\Sigma$,
$$
K(x,y) := K_x(y).
$$
\proclaim{Lemma \nextprocno}
\noindent For all $x\in\Bbb{H}^2$, the Green's function $K_x$ of $L$ over $\Bbb{H}^2$ with singularity at $x$ is
$$
K_x(y) = \frac{1}{2\pi} - \frac{1}{2\pi}\opArccoth(\opCosh(r(y)))\opCosh(r(y)),\eqnum{\nexteqnno[GreensFunctionA]}
$$
where $r(y)$ here denotes the distance in $\Bbb{H}^2$ from $y$ to $x$.
\endproclaim
\proclabel{GeneralGreensFunctionI}
\proof By uniqueness, $K_x$ is invariant under rotation about $x$ and is therefore a function of $r$ only. Since $L$ is given in polar coordinates of $\Bbb{H}^2$ about $x$ by
$$
Lu = u_{rr} + \opCoth(r)u_r + \frac{1}{\opSinh^2(r)}u_{\theta\theta} - 2u,
$$
$K_x$ is a solution of the equation
$$
F_{rr} + \opCoth(r)F_r - 2F = 0.
$$
We verify by inspection that the function $F(r):=\opCosh(r)$ is a solution and a second, linearly independent solution is then obtained using the Wronskian. The function $K_x$ is the unique linear combination of these two solutions which also satisfies
$$
\mlim_{r\rightarrow\infty}e^{2r}K_x(r) = -\frac{2}{3\pi},
$$
and
$$
\mlim_{r\rightarrow 0}r\frac{\partial}{\partial r}K_x(r) = \frac{1}{2\pi}.
$$
These properties imply that $K_x$ is an element of $L^1(\Bbb{H}^2)$ which solves \eqnref{GreensFunctionCondition}, and the result follows.\qed
\proclaim{Corollary \nextprocno}
\noindent The Green's kernel $K(x,y)$ of $L$ over $\Bbb{H}^2$ has the following properties.
\medskip
\myitem{(1)} For all $x\neq y$, $K(x,y)<0$;
\medskip
\myitem{(2)} for all $x\neq y$, $K(x,y)=K(y,x)$; and
\medskip
\myitem{(3)} for all $R>0$, there exists $C>0$ such that for $d(x,y)>R$,
$$
\left|K(x,y)\right| \leq Ce^{-2d(x,y)}.\eqnum{\nexteqnno[DecayRateOfPointGreensFunction]}
$$
\endproclaim
\proclabel{PropertiesOfGreensOperatorInHyperbolicSpace}
\noindent In particular,
\proclaim{Lemma \nextprocno}
\noindent Let $\Sigma$ be a hyperbolic surface. For all $R>0$, there exists $C>0$ with the property that
\medskip
\myitem{(1)} if $f$ is a twice-differentiable function such that both $f$ and $Lf$ are bounded; and
\medskip
\myitem{(2)} if $x\in\Sigma$ is such that $B_r(x)\minter\opSupp(Lf)=\emptyset$ for some $r\geq R$,
\medskip
\noindent then
$$
\left|f(x)\right|\leq Ce^{-r}\|Lf\|_{L^\infty}.\eqnum{\nexteqnno[DecayAwayFromSource]}
$$
\endproclaim
\proclabel{ExponentialDecayAwayFromSource}
\proof It suffices to consider the case where $\Sigma=\Bbb{H}^2$. Denote $g:=Lf$. Since $L$ has trivial kernel over $C^2_\opbdd(\Bbb{H}^2)$, we have
$$
f(x) = \int_{\Bbb{H}^2}K(x,y)g(y)\opdArea_y.
$$
However, in polar coordinates $(r,\theta)$ of $\Bbb{H}^2$ about $x$,
$$
\opdArea = \opSinh(r)drd\theta,
$$
and the result follows by Item $(3)$ of Corollary \procref{PropertiesOfGreensOperatorInHyperbolicSpace}.\qed
\newsubhead{Green's functions with geodesic singularities}[GreensFunctionsWithGeodesicSingularities]
Let $\Sigma$ be a complete hyperbolic surface without cusps and let $\gamma$ be a simple, closed geodesic in $\Sigma$. The {\sl Green's function} of $L$ over $\Sigma$ with singularity along $\gamma$ is defined to be the unique function $K_\gamma\in L^1(\Sigma)$ such that
$$
LK_\gamma\opdArea = \delta_\gamma\eqnum{\nexteqnno[GreensFunctionConditionII]}
$$
in the distributional sense, where $\delta_\gamma$ denotes the {\sl Dirac delta distribution} with singularity along $\gamma$, that is
$$
\langle\delta_\gamma,f\rangle := \int_\gamma f(x)dl_x.
$$
Using Fubini's Theorem, we readily verify
\proclaim{Lemma \nextprocno}
\noindent Let $\Sigma$ be a hyperbolic surface without cusps. For every simple, closed geodesic $\gamma$ in $\Sigma$, the Green's function $K_\gamma$ of $L$ over $\Sigma$ with singularity along $\gamma$ is given by
$$
K_\gamma(x) = \int_\gamma K(x,y)dl_y.
$$
\endproclaim
\proclabel{FirstFormulaForGreensFunction}
\noindent In order to derive more explicit estimates for $K_\gamma$, we consider the following special case. An {\sl hourglass} is defined to be a complete, connected hyperbolic surface containing a unique simple closed geodesic, that is, a complete hyperbolic annulus not conformal to the pointed disk.
\proclaim{Lemma \nextprocno}
\noindent Let $\Sigma$ be an hourglass with simple, closed geodesic $\gamma$. The Green's function $K_\gamma$ of $L$ over $\Sigma$ with singularity along $\gamma$ is given by
$$
K_\gamma(x) = -\frac{1}{\pi} + \frac{1}{\pi}\opArcCot(\opSinh(r(x)))\opSinh(r(x)),\eqnum{\nexteqnno[GreensFunctionB]}
$$
where $r(x)$ here denotes the distance in $\Sigma$ from $x$ to $\gamma$.
\endproclaim
\proclabel{FormulaForGreensFunctionInClepsydra}
\proof By uniqueness, $K_\gamma$ is invariant under rotation along and reflection about $\gamma$ and is therefore a function of $r$ only. In Fermi coordinates of $\Sigma$ about $\gamma$, $L$ is given by
$$
Lu = u_{rr} + \opTanh{r}u_r + \frac{1}{\opCosh^2(r)}u_{tt} - 2u,
$$
where $t$ here denotes a path-length parameter of $\gamma$. It follows that $K_\gamma$ is a solution of the equation
$$
F_{rr} + \opTanh(r)F_r - 2F = 0.
$$
We verify by inspection that the function $F(r):=\opSinh(r)$ is a solution and a second, linearly independent solution is then obtained using the Wronskian. The function $K_\gamma$ is the unique linear combination of these two functions which also satisfies
$$
\mlim_{r\rightarrow\infty}e^{2r}K_\gamma(r) = -\frac{2}{3\pi},
$$
and
$$
\mlim_{r\rightarrow 0^+}\frac{\partial}{\partial r}K_\gamma(r) - \mlim_{r\rightarrow 0^-}\frac{\partial}{\partial r}K_\gamma(r) = 1.
$$
These properties imply that $K_\gamma$ is an element of $L^1(\Sigma)$ which solves \eqnref{GreensFunctionConditionII}, and the result follows.\qed
\proclaim{Corollary \nextprocno}
\noindent Let $\Sigma$ be an hourglass with simple, closed geodesic $\gamma$. The Green's function $K_\gamma$ of $L$ over $\Sigma$ with singularity along $\gamma$ has the following properties.
\medskip
\myitem{(1)} For all $x$, $K_\gamma(x)<0$;
\medskip
\myitem{(2)} for all $x\in\gamma$, $K_\gamma(x)=-\frac{1}{\pi}$; and
\medskip
\myitem{(3)} there exists a constant $C>0$, which does not depend on $\Sigma$ such that, for all $x$,
$$
\left|K_\gamma(x)\right| \leq Ce^{-2d(\gamma,x)}.\eqnum{\nexteqnno[DecayRateOfGeodesicGreensFunction]}
$$
\endproclaim
\proclabel{PropertiesOfGreensFunctionGeodesic}
\proclaim{Lemma \nextprocno}
\noindent Let $\Sigma$ be a complete hyperbolic surface without cusps. Let $\gamma$ be a simple, closed geodesic in $\Sigma$. Let $\hat{\Sigma}$ be an hourglass with unique simple, closed geodesic $\hat{\gamma}$. Let $\pi:\hat{\Sigma}\rightarrow\Sigma$ be a local isometry whose restriction to $\hat{\gamma}$ defines an isometry onto $\gamma$. The Green's functions $K_\gamma$ and $K_{\hat{\gamma}}$ are related by
$$
K_\gamma(x) = \sum_{\pi(\hat{x})=x}K_{\hat{\gamma}}(\hat{x}).\eqnum{\nexteqnno[GeneralGeodesicGreensFunction]}
$$
\endproclaim
\proclabel{GeneralGreensFunctionGeodesic}
\proof It suffices to prove that the sum is locally uniformly absolutely convergent. However, for $x\in\Sigma$, the {\sl orbital counting function} of $x$ in $\hat{\Sigma}$ is defined by
$$
N(x)(R) := \# \pi^{-1}(\left\{x\right\})\minter B_R(\hat{\gamma}),
$$
where $B_R(\hat{\gamma})$ denotes the tubular neighbourhood of radius $R$ about $\hat{\gamma}$ in $\hat{\Sigma}$. By comparing areas, we obtain
$$
N(x)(R) \leq \frac{l\opSinh(R+r_\opinj)}{\opSinh^2(r_\opinj/2)},
$$
where $l$ here denotes the length of $\gamma$ and $r_{\opinj}$ denotes the injectivity radius of $\Sigma$ about $x$. Local uniform convergence follows from this and Item $(3)$ of Corollary \procref{PropertiesOfGreensFunctionGeodesic}. This completes the proof.\qed
\proclaim{Corollary \nextprocno}
\noindent Let $\Sigma$ be a complete hyperbolic surface without cusps. Let $\gamma$ be a simple, closed geodesic in $\Sigma$. The Green's function $K_\gamma$ of $L$ over $\Sigma$ with singularity along $\gamma$ has the following properties.
\medskip
\myitem{(1)} For all $x$, $K_\gamma(x)<0$; and
\medskip
\myitem{(2)} for all $x\in\gamma$, $K_\gamma(x)\leq -1/\pi$.
\endproclaim
\proclabel{PropertiesOfGeodesicGreensFunction}
\remark In fact, equality holds in the second relation at a single point if and only if $\Sigma$ is an hourglass.
\medskip
\noindent The following estimate will play a key role in the sequel.
\proclaim{Lemma \nextprocno}
\noindent For all $\epsilon>0$, there exists $C>0$ with the property that if $\Sigma$ is a complete hyperbolic surface without cusps, if $\gamma$ is a simple, closed geodesic in $\Sigma$, if $K_\gamma$ is the Green's function of $L$ over $\Sigma$ with singularity along $\gamma$, and if the injectivity radius of $\Sigma$ about every point of $\gamma$ is bounded below by $\epsilon$, then, for all $x\in\Sigma$,
$$
\left|K_\gamma(x)\right| \leq Ce^{-d(\gamma,x)}.\eqnum{\nexteqnno[DecayEstimateGeodesicGreensFunction]}
$$
\endproclaim
\proclabel{PropertiesOfGeodesicGreensFunctionII}
\proof Let $\pi:\Bbb{H}^2\rightarrow\Sigma$ be the canonical projection. Let $\hat{x}$ be an element of $\pi^{-1}(\left\{x\right\})$. Let $\Cal{G}$ be the set of complete geodesics in $\Bbb{H}^2$ which are preimages of $\gamma$ under $\pi$. The {\sl orbital counting function} $N(\gamma,x):[0,\infty[\rightarrow\Bbb{R}$ of $\gamma$ with respect to $x$ is given by
$$
N(\gamma,x)(R) := \#\left\{\hat{\gamma}\in\Cal{G}\ |\ \hat{\gamma}\minter B_R(\hat{x})\neq\emptyset\right\}.
$$
By comparing areas, we obtain
$$
N(\gamma,x)(R) \leq \frac{\opSinh^2((R+\epsilon)/2)}{\opSinh^2(\epsilon/2)}.
$$
By Lemmas \procref{FormulaForGreensFunctionInClepsydra} and \procref{GeneralGreensFunctionGeodesic},
$$
K_\gamma(x) = \sum_{\hat{\gamma}\in\Cal{G}}\bigg(-\frac{1}{\pi} + \frac{1}{\pi}\opArccot(\opSinh(d(\hat{\gamma},\hat{x})))\opSinh(d(\hat{\gamma},\hat{x}))\bigg).
$$
It follows by Item (3) of Corollary \procref{PropertiesOfGreensFunctionGeodesic} that, for some $C>0$,
$$\eqalign{
\left|K_\gamma(x)\right|
&\leq C\sum_{\hat{\gamma}\in\Cal{G}}e^{-2d(\hat{\gamma},\hat{x})}\cr
&\leq C\int_0^{e^{-2d(\gamma,x)}}N(\gamma,x)\bigg(-\frac{1}{2}\opLn(t)\bigg)dt\cr
&\leq \frac{Ce^\epsilon}{4\opSinh^2(\epsilon/2)}e^{-d(\gamma,x)},\cr}
$$
as desired.\qed
\newhead{Graftings and earthquakes}[GraftingsAndEarthquakes]
\newsubhead{Graftings and earthquakes}[GraftingsAndEarthquakes]
Let $\Sigma$ be a complete, marked hyperbolic surface without cusps. Let $g$ denote its metric. We recall the notation introduced in Section \subheadref{MainResults}. Given a simple, closed geodesic $\gamma$ in $\Sigma$, let $[\gamma]$ denote the free homotopy class in which it lies. Let $\l([\gamma],\Sigma)$ denote the infimal length with respect to $g$ amongst all curves in $[\gamma]$. Recall (see \cite{BallmanSchroederGromov}) that $[\gamma]$ contains no other geodesics and that $l([\gamma],\Sigma)$ is realised by $\gamma$.
\par
Recall that, given $t>0$, the {\sl grafting} of length $t$ of $\Sigma$ along $\gamma$, which we denote by $G(\gamma,\Sigma)(t)$, is defined to be the unique complete, marked hyperbolic surface obtained by cutting $\Sigma$ along $\gamma$, inserting a cylinder of length $t$, and multiplying the metric of the resulting surface by a suitable conformal factor. In this section, we show
\proclaim{Theorem \nextprocno}
\noindent Let $\Sigma$ be a complete, marked hyperbolic surface without cusps. Let $\gamma$ be a simple, closed geodesic in $\Sigma$, and for all $t\geq 0$, denote $\Sigma_t:=G(\gamma,\Sigma)(t)$. We have,
$$
\frac{\partial}{\partial t}\opLog(l([\gamma],\Sigma_t))\bigg|_{t=0} \leq -\frac{1}{\pi}.\eqnum{\nexteqnno[RatesOfChangeOfLengthI]}
$$
Furthermore, for all $\epsilon>0$, there exists $C>0$, which does not depend on $\Sigma$, such that, if the injectivity radius of $\Sigma$ is bounded below by $\epsilon>0$, then, for every other simple, closed geodesic $\gamma'$ not intersecting $\gamma$,
$$
\bigg|\frac{\partial}{\partial t}\opLog(l([\gamma'],\Sigma_t))\bigg|_{t=0}\bigg| \leq Ce^{-d(\gamma,\gamma')}.\eqnum{\nexteqnno[RatesOfChangeOfLengthII]}
$$
\endproclaim
\proclabel{RatesOfChangeOfLength}
Theorem \procref{RatesOfChangeOfLength} will be proven at the end of the following section. We first show how graftings are also constructed as smooth non-conformal perturbations of $g$. For completeness, we also study earthquakes, as their analysis is almost identical, and this will allow us to recover the result \cite{Wolpert} of Wolpert. As outlined in Section \subheadref{MainResults}, given $t\in\Bbb{R}$, the (right) {\sl earthquake} of length $t$ of $\Sigma$ along $\gamma$ of length $t$, denoted by $E(\gamma,\Sigma)(t)$, is defined to be the unique marked hyperbolic surface obtained by cutting $\Sigma$ along $\gamma$ and rotating the right hand side of $\gamma$ by a distance $t$ in the positive direction of this geodesic. Since the resulting metric is automatically hyperbolic, there is no need to multiply by a conformal factor in this case. This definition is independent of the orientation of $\gamma$ chosen, but depends on the orientation of the ambient surface. The left earthquake of $(\Sigma,g)$ along $\gamma$ of length $t$ is defined in a similar manner and is readily shown to be equal to $E(\gamma,\Sigma)(-t)$. Reversing the orientation of $\Sigma$ interchanges left and right earthquakes.
\par
Let $R>0$ be such that the tubular neighbourhood $B_R(\gamma)$ of radius $R$ about $\gamma$ is isometric to $\gamma\times]-R,R[$ furnished with the twisted product metric
$$
\opCosh^2(r)dt^2 + dr^2.\eqnum{\nexteqnno[TwistedMetric]}
$$
For $\phi,\psi\in C_0^\infty(]-R,R[)$, define the metric $g_{\phi,\psi}$ over this annulus by
$$
g_{\phi,\psi} := (\opCosh(r)dt - \phi(r)dr)^2 + e^{2\psi(r)}dr^2,\eqnum{\nexteqnno[PerturbedMetric]}
$$
and extend $g_{\phi,\psi}$ to a smooth metric over the whole of $\Sigma$ by setting it equal to $g$ over the complement of this annulus.
\proclaim{Lemma \nextprocno}
\myitem{(1)} For all $\psi$, $(\Sigma,g_{0,\psi})$ is conformally equivalent to $G(\gamma,\Sigma)(t)$, where
$$
t = \int_{-R}^R\frac{e^{\psi(r)} - 1}{\opCosh(r)}dr.\eqnum{\nexteqnno[SizeOfGrafting]}
$$
\myitem{(2)} For all $\phi$, $(\Sigma,g_{\phi,0})$ is conformally equivalent to $E(\gamma,\Sigma)(t)$, where
$$
t = \int_{-R}^R\frac{\phi(s)}{\opCosh(s)}ds.\eqnum{\nexteqnno[SizeOfEarthquake]}
$$
\endproclaim
\proclabel{ModulesOfSmoothTransformations}
\remark The proof of Lemma \procref{ModulesOfSmoothTransformations} uses the concept of conformal module (see \cite{Ahlfors}). Recall that every annulus $A$ of hyperbolic type is conformally equivalent to $S^1\times]0,M[$ for a unique $M\in]0,\infty]$. For our purposes, the {\sl conformal module} of $A$ is defined to be equal to this number.
\medskip
\proof Consider first the map
$$
\Psi:N_R(\gamma)\rightarrow S^1\times]0,M[;(t,r)\mapsto\frac{2\pi}{l}\bigg(t,\int_{-R}^r\frac{1}{\opCosh(s)}dr\bigg),
$$
where $l$ here denotes the length of $\gamma$ and $M$ is to be determined. We verify that $\Phi$ is conformal with respect to the metric $g$ over $B_R(\gamma)$ and the product metric over $S^1\times ]0,M[$. Thus, if we denote by $M(B_R(\gamma),g)$ the conformal module of this cylinder, then
$$
M(B_R(\gamma),g) = \frac{2\pi}{l}\int_{-R}^R\frac{1}{\opCosh(s)}ds.
$$
The conformal module of the cylinder $(B_R(\gamma),g_{0,\psi})$ is likewise given by
$$
M(B_R(\gamma),g_{0,\psi}) = \frac{2\pi}{l}\int_{-R}^R\frac{e^{\psi(s)}}{\opCosh(s)}ds.
$$
On the other hand, for all $t$, the conformal module of the cylinder $(\gamma\times]0,t[,dt^2+dr^2)$ is
$$
M(\gamma\times]0,t[,dt^2+dr^2) = \frac{2\pi t}{l}.
$$
Item $(1)$ now follows upon comparing these three conformal modules. Finally, for $\phi\in C_0^\infty(]-R,R[)$, define
$$
\Phi_\phi:B_R(\gamma)\rightarrow B_R(\gamma);(t,r)\mapsto\bigg(t-\int_{-R}^R\frac{\phi(s)}{\opCosh(s)}dr,r\bigg).
$$
Since $\Phi_\phi$ defines an isometry from $(B_R(\gamma),g_{\phi,0})$ to $(B_R(\gamma),g)$, Item $(2)$ follows, and this completes the proof.\qed
\newsubhead{First order variations of curvature}[FirstOrderVariationsOfCurvatureAndLength]
Let $\Sigma$ be a complete, marked hyperbolic surface without cusps and let $g$ denote its metric. Let $\gamma$ be a simple, closed geodesic in $\Sigma$. For suitable $R>0$ and for $\phi,\psi\in C_0^2(]-R,R[)$, let $g_{\phi,\psi}$ denote the complete metric obtained by perturbing $g$ in a neighbourhood of $\gamma$ as in Section \subheadref{GraftingsAndEarthquakes}.
\par
For $\phi,\psi\in C_0^\infty(]-R,R[)$ and for $f\in C^2(\Sigma)$, denote
$$
g_{\phi,\psi,f} := e^{2f}g_{\phi,\psi},\eqnum{\nexteqnno[ThreeParameterPerturbedMetric]}
$$
and let $\kappa_{\phi,\psi,f}$ denote the curvature function of this metric. The operator $\kappa$ defines a smooth functional from a neighbourhood of $(0,0,0)$ in $C_0^{2,\alpha}(]-R,R[)^2\times C^{2,\alpha}(\Sigma)$ into $C^{0,\alpha}(\Sigma)$. Furthermore, for all $f$, we have (c.f. \cite{ChowKnopf}),
$$
\kappa_{0,0,f} = -e^{-2f}\big(\Delta f + 1\big).\eqnum{\nexteqnno[ConformalChangeOfCurvature]}
$$
In particular, considered as a smooth map between open subsets of Banach spaces, the partial derivative of $\kappa$ at $(0,0,0)$ with respect to the third component is given by
$$
D_3\kappa_{0,0,0}f = -Lf,\eqnum{\nexteqnno[DerivativeOfCurvatureWithRespectToThirdComponent]}
$$
where $L$ is the operator introduced in Section \subheadref{GreensFunctionsWithPointSingularities}. The existence of a Green's kernel ensures that $L$ defines a linear isomorphism from $C^{2,\alpha}(\Sigma)$ into $C^{0,\alpha}(\Sigma)$ (c.f. \cite{GilbTrud}). The Implicit Function Theorem then yields
\proclaim{Lemma \nextprocno}
\noindent For all $\phi,\psi\in C_0^{2,\alpha}(]-R,R[)$, there exists $\epsilon>0$ and a smooth function $f(\phi,\psi):]-\epsilon,\epsilon[\rightarrow C^{2,\alpha}(\Sigma)$ such that, for all $s$,
$$
\kappa_{s\phi,s\psi,f(\phi,\psi)(s)} = -1.
$$
\endproclaim
Let $\gamma'$ be another simple, closed geodesic in $\Sigma$. For all $(\phi,\psi,f)$, let $l_{\phi,\psi,f}([\gamma'])$ denote the infimal length of curves in $[\gamma']$ with respect to the metric $g_{\phi,\psi,f}$.
\proclaim{Lemma \nextprocno}
\noindent For all $(\phi,\psi)\in C_0^{2,\alpha}(]-R,R[)$,
$$\eqalign{
\left.\frac{\partial}{\partial s}l_{s\phi,s\psi,f(\phi,\psi)(s)}([\gamma'])\right|_{s=0} &
=-\int_{\gamma'}(\phi\circ r)\opCosh^2(r)\frac{dt}{dl'}\frac{dr}{dl'}dl'+\int_{\gamma'}(\psi\circ r)\bigg(\frac{dr}{dl'}\bigg)^2 dl'\cr
&\qquad\qquad+\int_{\gamma'}\left.\frac{\partial}{\partial s}f(\phi,\psi)(s)\right|_{s=0}dl',\cr}
\eqnum{\nexteqnno[ElementaryFormulaForDerivativeOfLength]}
$$
where $(t,r)$ are the coordinates of $B_R(\gamma)$ given in Section \subheadref{GraftingsAndEarthquakes}.
\endproclaim
\proclabel{ElementaryFormulaForDerivativeOfLength}
\remark Since $(\psi\circ r)$ and $(\phi\circ r)$ are supported in $B_R(\gamma)$, the integrands of the first two terms on the right-hand side of \eqnref{ElementaryFormulaForDerivativeOfLength} are non-trivial only along those segments of $\gamma'$ which lie inside this tubular neighbourhood. In particular, since $r$ is smooth over these segments, these terms are indeed well-defined.
\medskip
\proof For all $(\phi,\psi,f)$ and for all $\eta\in C^{2,\alpha}(\gamma',]-\delta,\delta[)$, let $l'_{\phi,\psi,f,\eta}$ denote the length of the graph of $\eta$ over $\gamma'$ with respect to the metric $g_{\phi,\psi,f}$. This is a smooth functional whose partial derivatives with respect to the first three components at zero are
$$\eqalign{
D_1l_{0,0,0,0}\phi &= -\int_{\gamma'}(\phi\circ r)\opCosh^2(r)\frac{dt}{dl'}\frac{dr}{dl'}dl',\cr
D_2l_{0,0,0,0}\psi &= \int_{\gamma'}(\psi\circ r)\bigg(\frac{dr}{dl'}\bigg)^2dl',\ \text{and}\cr
D_3l_{0,0,0,0}f &= \int_{\gamma'} fdl'.}
$$
Since $\gamma'$ is a critical point of the length functional for $g$, its partial derivative with respect to the fourth component at zero vanishes. Since geodesics in hyperbolic surfaces are stable under small perturbations, upon decreasing $\epsilon$ if necessary, there exists a smooth function $\eta:]-\epsilon,\epsilon[\rightarrow C^{2,\alpha}(\gamma',]-\delta,\delta[)$ such that, for all $s$, the graph of $\eta(s)$ is the unique geodesic in $[\gamma']$ which realises the infimal length with respect to the metric $g_{s\phi,s\psi,f(\phi,\psi)(s)}$ amongst curves in this class. The result now follows by the chain rule.\qed
\proclaim{Lemma \nextprocno}
\noindent Over $B_R(\gamma)$, the curvature of $g_{\phi,\psi,0}$ is given by
$$
\kappa_{\phi,\psi,0} = e^{-2\psi(r)}\psi'(r)\opTanh(r) - e^{-2\psi(r)}.\eqnum{\nexteqnno[PerturbedCurvatureFormula]}
$$
\endproclaim
\proclabel{PerturbedCurvatureFormula}
\proof Consider first a general metric $h$ over $\Sigma$. Let $(e_1,e_2)$ be a local orthonormal moving frame of $h$. Recall (c.f. \cite{Chavel} and \cite{Spivak}) that, for each $i$ and for every vector field $\xi$,
$$
\nabla_\xi^he_i = -\alpha_h(\xi)J^he_i,
$$
where $J^h$ denotes the complex structure of $h$, $\nabla^h$ denotes its Levi-Civita covariant derivative and the {\sl connection form} $\alpha_h$ is given by
$$
\alpha_h(\xi) := h([e_1,e_2],\xi).\eqnum{\nexteqnno[ConnectionFormFormula]}
$$
Recall that the curvature of $h$ is then given by
$$
\kappa_h = d\alpha_h(e_1,e_2).\eqnum{\nexteqnno[MovingFrameCurvature]}
$$
\par
Consider now the local orthonormal moving frame given by
$$\eqalign{
e_1 &:= \frac{1}{\opCosh(r)}\partial_t,\hfill\cr
e_2 &:= e^{-\psi(r)}\partial_r + \frac{\phi(r)e^{-\psi(r)}}{\opCosh(r)}\partial_t.\hfill\cr}
$$
We compute
$$
[e_1,e_2] = e^{-\psi(r)}\opTanh(r)e_1,
$$
and the result now follows by \eqnref{MovingFrameCurvature}.\qed
\medskip
Consider now sequences $(R_m)$, $(\phi_m)$ and $(\psi_m)$ such that,
\medskip
\myitem{(1)} $(R_m)\downarrow 0$ as $m$ tends to infinity;
\medskip
\myitem{(2)} $\phi_m,\psi_m\in C_0^\infty(]-R_m,R_m[)$;
\medskip
\myitem{(3)} $\int_{-R_m}^{R_m}\opCosh(s)\psi_m(s)ds = 1$; and
\medskip
\myitem{(4)} $(\|\psi_m\|_{L^1})$ is uniformly bounded independent of $m$.
\medskip
\noindent For all $m$, denote
$$
\tilde{\kappa}_m := \left.\frac{\partial}{\partial s}\kappa_{s\phi_m,s\psi_m,0}\right|_{s=0}.\eqnum{\nexteqnno[DerivativeOfCurvature]}
$$
\proclaim{Lemma \nextprocno}
\noindent For any Lipschitz function $f$,
$$
\mlim_{m\rightarrow\infty}\int_{\Sigma}f\tilde{\kappa}_m\opdArea = \int_\gamma fdl.\eqnum{\nexteqnno[DistributionalLimitOfDerivativeOfCurvature]}
$$
\endproclaim
\proclabel{GaussBonnet}
\remark In other words, $\tilde{\kappa}_m\opdArea$ converges towards $\delta_\gamma$ in the distributional sense as $m$ tends to infinity.
\medskip
\proof Indeed, for all $m$ and for all sufficiently small $s$, denote $g_{m,s}:=g_{s\phi_m,s\psi_m,0}$ and $\kappa_{m,s}:=\kappa_{s\phi_m,s\psi_m,0}$. By \eqnref{PerturbedCurvatureFormula},
$$
\kappa_{m,s}(r)= s\psi_m'(r)\opTanh(r)e^{-2s\psi_m(r)}-e^{-2s\psi_m(r)}.
$$
Differentiating this relation with respect to $s$ at $s=0$ yields
$$
\tilde{\kappa}_m(r) = 2\psi_m(r) + \psi'_m(r)\opTanh(r).
$$
Bearing in mind Property $(3)$ of the sequence $(\psi_m)$, upon integrating by parts we obtain, for Lipschitz $f$,
$$\eqalign{
\int_{-R_m}^{R_m}f(t,r)\tilde{\kappa}_m(r)\opCosh(r)dr
&=2\int_{-R_m}^{R_m}f(t,r)\psi_m(r)\opCosh(r)dr\cr
&\qquad\qquad + \int_{-R_m}^{R_m}f(t,r)\psi'_m(r)\opSinh(r)dr\cr
&=f(t,0) + \int_{-R_m}^{R_m}\big(f(t,r) - f(t,0)\big)\psi_m(r)\opCosh(r)dr\cr
&\qquad\qquad - \int_{-R_m}^{R_m}\psi_m(r)\opSinh(r)f_r(t,r)dr,\cr}
$$
where the derivative of $f$ with respect to $r$ is taken in the distributional sense. However, for all $m$,
$$\eqalign{
\bigg|\int_{-R_m}^{R_m}\big(f(t,r) - f(t,0)\big)\psi_m(r)\opCosh(r)dr\bigg| &\leq 2R_m\opCosh(R_m)\|\psi_m\|_{L^1}[f]_1,\ \text{and}\cr
\bigg|\int_{-R_m}^{R_m}\psi_m(r)\opSinh(r)f'(r)dr\bigg| &\leq \opSinh(R_m)\|\psi_m\|_{L^1}[f]_1,\cr}
$$
where $[f]_1$ here denotes the Lipschitz seminorm of $f$. The result follows upon integrating with respect to $t$ and letting $m$ tend to infinity.\qed
\proclaim{Theorem \nextprocno}
\noindent For every pair $(\gamma,\gamma')$ of simple, closed geodesics in $\Sigma$,
$$\eqalignno{
\left.\frac{\partial}{\partial t}l([\gamma'],G(\gamma,\Sigma)(t))\right|_{t=0}
&=\sum_{x\in\gamma\minter\gamma'}\opSin(\theta_x) + \int_\gamma\int_{\gamma'}K(p,q)dl_pdl_q,\ \text{and}&\nexteqnno[GraftingVariationFormula]\cr
\left.\frac{\partial}{\partial t}l([\gamma'],E(\gamma,\Sigma)(t))\right|_{t=0}
&=\sum_{x\in\gamma\minter\gamma'}\opCos(\theta_x),&\nexteqnno[EarthquakeVariationFormula]\cr}
$$
where $K(x,y)$ is the Green's kernel of $L$ over $\Sigma$ and, for all $x\in\gamma\minter\gamma'$, $\theta_x$ denotes the angle that $\gamma$ makes with $\gamma'$ at the point $x$.
\endproclaim
\remark Equation \eqnref{EarthquakeVariationFormula} was first determined by Wolpert in \cite{Wolpert}.
\medskip
\remark By convention, the sums on the right-hand sides of \eqnref{GraftingVariationFormula} and \eqnref{EarthquakeVariationFormula} are taken to be zero when $\gamma=\gamma'$.
\medskip
\remark If $\gamma\neq\gamma'$, then these two geodesics intersect transversally, so that the above sums are finite and thus indeed well-defined.
\medskip
\proof Let $\psi\in C_0^\infty(]-R,R[)$ be a positive function such that
$$
\int_{-R}^R\frac{\psi(r)}{\opCosh(r)}dr = 1.
$$
By Lemma \procref{ModulesOfSmoothTransformations}, up to first order in $t$, $(\Sigma,g_{0,t\psi})$ is conformally equivalent to $G(\gamma,\Sigma)(t)$. For all $t$, let $f_t:=f(t)\in C^{2,\alpha}(\Sigma)$ be the unique function such that $g_{0,t\psi,f(t)}$ is hyperbolic. Then, by Lemma \procref{ElementaryFormulaForDerivativeOfLength},
$$
\frac{\partial}{\partial t}l([\gamma'],G(\gamma,\Sigma)(t))\bigg|_{t=0}=\int_{\gamma'}(\psi\circ r)(x)\bigg(\frac{dr}{dl'}\bigg)^2dl_x' + \int_{\gamma'}\frac{\partial}{\partial t}f(x)\bigg|_{t=0}dl_x'.
$$
Letting $R$ tend to zero in the first integral yields
$$
\mlim_{R\rightarrow 0}\int_{\gamma'}\psi(r)\bigg(\frac{dr}{dl'}\bigg)^2dl_x' =
\sum_{x\in\gamma\minter\gamma'}\frac{dr}{dl'}(x) =
\sum_{x\in\gamma\minter\gamma'}\opSin(\theta_x).
$$
On the other hand, for all $x$,
$$
\left.\frac{\partial}{\partial t}\kappa_{0,t\psi,0}(x)\right|_{t=0}-\left.L\frac{\partial}{\partial t}f_t(x)\right|_{t=0} =
\left.\frac{\partial}{\partial t}\kappa_{0,t\psi,f(t)}\right|_{t=0} = 0,
$$
so that,
$$
\left.\frac{\partial}{\partial t}f_t(x)\right|_{t=0} = \int_\Sigma K(x,y)\tilde{\kappa}_{\psi}(y)dArea_y,
$$
where
$$
\tilde{\kappa}_\psi = \left.\frac{\partial}{\partial t}\kappa_{0,t\psi,0}(y)\right|_{t=0}.
$$
Thus, bearing in mind Fubini's theorem, Item $(2)$ of Corollary \procref{PropertiesOfGreensOperatorInHyperbolicSpace} and Lemma \procref{FirstFormulaForGreensFunction},
$$\eqalign{
\int_{\gamma'}\frac{\partial}{\partial t}f_t(x)\bigg|_{t=0}dl'_x
&=\int_{\gamma'}\int_\Sigma K(x,y)\tilde{\kappa}_\psi(y)dArea_y dl'_x\cr
&=\int_\Sigma\int_{\gamma'}K(y,x)\tilde{\kappa}_\psi(y)dl'_x dArea_y\cr
&=\int_\Sigma K_{\gamma'}(y)\tilde{\kappa}_\psi(y)dArea_y.\cr}
$$
Since $K_{\gamma'}$ is Lipschitz, it follows by Lemma \procref{FirstFormulaForGreensFunction} again and Lemma \procref{GaussBonnet} that
$$
\mlim_{R\rightarrow 0}\int_{\gamma'}\left.\frac{\partial}{\partial t}f_t(x)\right|_{t=0}dl_x = \int_\gamma K_{\gamma'}(x)dl_x = \int_{\gamma}\int_{\gamma'}K(x,y)dl_y'dl_x,
$$
and the first relation follows. The second relation is proven in an analogous manner, and this completes the proof.\qed
\medskip
{\bf\noindent Proof of Theorem \procref{RatesOfChangeOfLength}:} When $\gamma'=\gamma$ or when $\gamma'$ is disjoint from $\gamma$, the first term on the right hand side of \eqnref{GraftingVariationFormula} vanishes, so that, bearing in mind Lemma \procref{FirstFormulaForGreensFunction},
$$
\left.\frac{\partial}{\partial t}l([\gamma'],G(\gamma,\Sigma)(t))\right|_{t=0}
=\int_{\gamma'}\int_\gamma K(x,y)dl_xdl_y
=\int_{\gamma'}K_\gamma(x)dl_x.
$$
Thus, by Item $(2)$ of Corollary \procref{PropertiesOfGeodesicGreensFunction},
$$
\left.\frac{\partial}{\partial t}l([\gamma],G(\gamma,\Sigma)(t))\right|_{t=0}
\leq -\int_\gamma\frac{1}{\pi}dl_x
=-\frac{l([\gamma],\Sigma)}{\pi},
$$
and the first result follows. Likewise, with $C$ as in Item $(3)$ of Corollary \procref{PropertiesOfGeodesicGreensFunction},
$$
\left|\left.\frac{\partial}{\partial t}l([\gamma'],G(\gamma,\Sigma)(t))\right|_{t=0}\right|
\leq C\int_{\gamma'}e^{-d(\gamma,\gamma')}dl_x
=Cl([\gamma'],\Sigma)e^{-d(\gamma,\gamma')},
$$
and the second result follows. This completes the proof.\qed
\newhead{Laminations and hyperbolic perturbations}[Laminations]
\newsubhead{Riemannian laminations}[Laminations]
We now recall basic definitions and results of the theory of laminations, referring the reader to \cite{CandelConlon} for a thorough introduction. Let $X$ be a topological space. A $d$-dimensional {\sl laminated chart} of $X$ is defined to be a pair $(U,T,\Phi)$ where $U$ is an open subset of $X$, $T$ is a topological space and $\Phi:U\rightarrow]-1,1[^d\times T$ is a homeomorphism. Given two laminated charts $(U_i,T_i,\Phi_i)_{i\in\left\{1,2\right\}}$, the transition map $\alpha_{21}:\Phi_1(U_1\minter U_2)\rightarrow\Phi_2(U_1\minter U_2)$ is defined by
$$
\alpha_{21} := \Phi_2\circ\Phi_1^{-1}.
$$
This map is said to be of class $C^\infty_l$ whenever every point of $\Phi_1(U_1\minter U_2)$ has a neighbourhood $\Omega$ of the form
$$
\Omega := ]x_1-\epsilon,x_1+\epsilon[\times...\times]x_d-\epsilon,x_d+\epsilon[\times S
$$
over which
$$
\alpha_{21}(x,t) = (\phi_{21}(x,t),\tau_{21}(t)),
$$
where
\medskip
\myitem{(1)} $\tau_{21}$ is a homeomorphism onto its image;
\medskip
\myitem{(2)} for all $t$, $\phi_{21}(\cdot,t)$ is a smooth diffeomorphism onto its image; and
\medskip
\myitem{(3)} $\phi_{21}(\cdot,t)$ varies continuously with $t$ in the $C^\infty_\oploc$ sense.
\medskip
\noindent A {\sl lamination} is defined to be a separable, metrizable space $X$ furnished with an atlas $\Cal{A}$ of laminated charts all of whose transition maps are of class $C^\infty_l$.
\par
Given a lamination $X$, a laminated chart $(U,T,\Phi)$ of $X$ and a point $t$ of $T$, we call the set $\Phi^{-1}(]-1,1[^d\times\left\{t\right\})$ a {\sl plaque} of the chart. Plaques glue together to yield a partition of $X$ into smooth, $d$-dimensional manifolds, called {\sl leaves} of the lamination. For all $x\in X$, the leaf passing through $x$ will be denoted by $\Sigma_x$. Given two laminations $X$ and $Y$ and $k\in\Bbb{N}\munion\left\{\infty\right\}$, $C^k_l(X,Y)$ is defined to be the space of all continuous functions $f:X\rightarrow Y$ which restrict to $C^k$ functions from leaves to leaves and which, in addition vary continuously in the $C^k_\oploc$ sense as the plaques vary in any given laminated chart. Functions in $C^\infty_l(X,Y)$ are said to be {\sl leafwise smooth}. Leafwise smooth maps are the morphisms of the category of laminations.
\par
Since every finite-dimensional manifold naturally carries the structure of a lamination consisting of a single leaf, the theory of laminations may be viewed as an extension of the theory of manifolds. Standard constructions of manifold theory then carry over to the theory of laminations with appropriate modifications. In particular, vector bundles over laminations are defined in the natural manner and, given a lamination $X$, the {\sl tangent bundle} $TX$ of $X$ is defined to be the vector bundle whose fibre at the point $x\in X$ is the tangent space to the leaf $\Sigma_x$ at this point. The cotangent bundle and other tensor bundles over $X$ are likewise defined in the natural manner.
\par
Given a lamination $X$ and $k\in\Bbb{N}\munion\left\{\infty\right\}$, let $C^k_l(X)$ denote the space of $C^k_l$ functions from $X$ into the trivial lamination $\Bbb{R}$. We recall
\proclaim{Theorem \nextprocno, Candel \cite{CandelI}}
\noindent For every open cover $(U_i)_{i\in I}$ of $X$, there exists a locally finite, leafwise smooth partition of unity $(\chi_j)_{j\in J}$ of $X$ subordinate to this cover.
\endproclaim
\proclabel{PartitionOfUnity}
\noindent When, in addition, $X$ is compact, leafwise smooth partitions of unity serve to furnish $C^k_l(X)$ with a canonical Banach space structure as follows. First, let $(U,T,\Phi)$ be a laminated chart of $X$ and, for $f\in C^k_l(]-1,1[^d\times T)$, define
$$
\|f\|_{C^k_l} := \msup_{t\in T}\|f(\cdot,t)\|_{C^k}.
$$
Next, let $(U_\alpha,T_\alpha,\Phi_\alpha)_{\alpha\in A}$ be a finite atlas of $X$ by laminated charts, let $(\chi_\alpha)_{\alpha\in A}$ be a leafwise smooth partition of unity of $X$ subordinate to this atlas and, for $f\in C^k_l(X)$, define
$$
\|f\|_{C^k_l} := \sum_{\alpha\in A}\|(f\chi_\alpha)\circ\Phi_\alpha^{-1}\|_{C^k_l}.\eqnum{\nexteqnno[LeafwiseBanachNorm]}
$$
Up to uniform equivalence, the norm \eqnref{LeafwiseBanachNorm} is independent of the atlas and partition of unity chosen, and thus defines a canonical Banach space structure over $C^k_l(X)$. For all $(k,\alpha)$, the H\"older norm $\|\cdot\|_{C^{k,\alpha}_l}$ is defined in a similar manner, and the H\"older space $C^{k,\alpha}_l(X)$ is defined to be the Banach space of all functions $f\in C^k_l(X)$ such that
$$
\|f\|_{C^{k,\alpha}_l(X)}<\infty.
$$
Finally, given any vector bundle $EX$ over $X$, for all $(k,\alpha)$, the Banach space of $C^{k,\alpha}_l$ sections of $EX$ over $X$ is defined in an analogous manner and will be denoted by $\Gamma^{k,\alpha}_l(EX)$.
\par
A {\sl leafwise metric} over $X$ is defined to be a leafwise smooth, positive-definite section of the bundle $\opSymm(TX)$ of symmetric, bilinear forms over $X$. In view of Theorem \procref{PartitionOfUnity}, leavewise metrics always exist and are constructed in the same way as in the classical theory of riemannian manifolds. Given a leafwise metric $g$, the pair $(X,g)$ will be called a {\sl riemannian lamination}. For all $x$, the restriction of $g$ to the leaf $\Sigma_x$ will be denoted by $g_x$.
\par
A {\sl hyperbolic surface lamination} is a riemannian lamination all of whose leaves are complete hyperbolic surfaces. In \cite{CandelI}, Candel characterises hyperbolic surface laminations in terms of the conformal classes of its leaves. Recall first that two leafwise metrics $g$ and $g'$ are said to be {\sl conformally equivalent} whenever there exists a leafwise smooth function $u$ such that
$$
g' = e^{2u}g.
$$
\proclaim{Theorem \nextprocno, Candel \cite{CandelI}}
\noindent Let $(X,g)$ be a compact riemannian surface lamination. If every leaf of $(X,g)$ is of hyperbolic type, then there exists a unique leafwise smooth metric $\tilde{g}$ in the conformal class of $g$ with respect to which every leaf is complete and hyperbolic.
\endproclaim
\proclabel{Candel}
\noindent A closer examination of Candel's proof yields
\proclaim{Lemma \nextprocno, Candel \cite{CandelI}}
\noindent Let $\Cal{G}^{k,\alpha}_l(X)$ denote the space of negatively curved, leafwise smooth metrics of class $C^{k,\alpha}_l$ over $X$. Let $\Cal{H}:\Cal{G}^{k,\alpha}_l(X)\rightarrow\Cal{G}^{k,\alpha}_l(X)$ be such that, for all $g$, $\Cal{H}(g)$ is the unique leafwise smooth hyperbolic metric given by Theorem \procref{Candel}. Then $\Cal{H}$ is smooth as a map between Banach manifolds.
\endproclaim
\proclabel{CandelSmoothness}
\newsubhead{Elementary differential topology of laminations}[TheDifferentialTopologyOfLaminations]
Let $X$ be a compact riemannian lamination. We gather here various elementary properties of $X$ which will be required in the sequel. First, the {\sl leafwise distance} over $X$ is defined by
$$
d_l(x,y) := \minf_{\gamma}l_g(\gamma),\eqnum{\nexteqnno[DefinitionOfLeafwiseDistance]}
$$
where $\gamma$ varies over all tangential curves from $x$ to $y$ and $l_g(\gamma)$ denotes the length of $\gamma$ with respect to $g$. In particular, the leafwise distance between two points is finite if and only if they both lie on the same leaf. This distance function defines the {\sl leafwise topology} of $X$, which is the smallest topology containing all open subsets of leaves of $X$. A sequence $(x_m)$ of points in $X$ converges to the point $x_\infty$ in this topology if and only if $(x_m)$ converges to $x_\infty$ in the ambient topology of $X$ and, for sufficiently large $m$, $x_m$ also lies in the same leaf as $x_\infty$. A subset $Y$ of $X$ is compact with respect to this topology if and only if it consists of a finite union of compact subsets of leaves.
\proclaim{Lemma \nextprocno}
\noindent Let $(U,T,\Phi)$ be a laminated chart of $X$. Then $T$ is separable, metrizable and locally compact.
\endproclaim
\proclabel{SepMetComp}
\proof Indeed, being a subset of a separable, metrizable space, $U$ is also separable and metrizable. Since $T=\left\{0\right\}\times T$ is homeomorphic to a subset of $U$, separability and metrizability follow. It remains to prove local compactness. However, choose $t\in T$. Denote $x:=\Phi^{-1}(\left\{0\right\})$. Since $x$ is metrizable, there exists a neighbourhood $V$ of $x$ in $X$ such that $\overline{V}\subseteq U$. Since $X$ is compact, so too is $\overline{V}$. Since $\Phi$ is a homeomorphism, $\Phi(V)$ is open and $\overline{\Phi{V}}=\Phi(\overline{V})$ is compact. Let $W$ be a neighbourhood of $t$ in $T$ such that $\left\{0\right\}\times W\subseteq\Phi(V)$. Then $\overline{W}$ is compact, and the result follows.\qed
\proclaim{Lemma \nextprocno}
\noindent For all $x\in X$, there exists $C>0$ and a laminated chart of $X$ containing $x$ whose plaques have volume bounded below by $1/C$ and diameter bounded above by $C$.
\endproclaim
\proclabel{BoundedAreaAndDiameter}
\remark By compactness of $X$, $B$ may even be chosen independent of $x$.
\medskip
\proof Let $(U,T,\Phi)$ be a laminated chart of $X$ about $x$. Let $t\in T$ be such that $\Phi(x)$ is contained in the plaque $]-1,1[^d\times\left\{t\right\}$. Let $V$ be a neighbourhood of $t$ in $T$ with compact closure.\qed
\proclaim{Lemma \nextprocno}
\noindent Let $Y\subseteq X$ be compact in the leafwise topology. For all $x\in X$, there exists a laminated chart $(U,T,\Phi)$ of $X$ about $x$ such that $U\minter Y$ is contained in, at most, a single plaque.
\endproclaim
\proclabel{SinglePlaque}
\proof We suppose that $X$ is furnished with a leafwise metric $g$. Let $(U,T,\Phi)$ be a laminated chart of $X$ about $x$ whose plaques have volume bounded below by $1/C$, say, and diameter bounded above by $C$, say. Since $Y$ is compact in the leafwise topology, there exist $y_1,...,y_m\in Y$ and $R>0$ such that
$$
Y\subseteq\munion_{i=1}^m B_R(y_i),
$$
where, for each $i$, $B_R(y_i)$ here denotes the ball of radius $R$ about $y_i$ in $X$ with respect to the leafwise metric. If $V$ is the volume of the set
$$
\munion_{i=1}^m B_{R+C}(y_i),
$$
then $(U,T,\Phi)$ contains at most $\lfloor CV\rfloor$ plaques which intersect $Y$ non-trivially where, for all $\lambda$, $\lfloor\lambda\rfloor$ denotes the greatest integer less than or equal to $\lambda$. The result now follows upon reducing $T$ if necessary.\qed
\proclaim{Lemma \nextprocno}
\noindent Let $Y\subseteq X$ be compact in the leafwise topology. Let $Z$ be the (finite) union of all leaves which intersect $Y$ non-trivially. There exists a neighbourhood $\Omega$ of $Y$ in $X$ and a leafwise smooth function $\pi:\Omega\rightarrow Z$ such that, for all $y\in Y$, $\pi(y)=y$.
\endproclaim
\proclabel{TubularNeighbourhood}
\proof It suffices to consider the case where $Y$ is contained in a single leaf $\Sigma$, say. By compactness, $Y$ is covered by a finite family $(U_i,T_i,\Phi_i)_{1\leq i\leq m}$ of laminated charts. By Lemma \procref{SinglePlaque}, we may suppose that, for all $i$,
$$
\Phi_i(Y\minter U_i) \subseteq ]-1,1[^d\times\left\{t_i\right\},
$$
for some $t_i\in T_i$. For all $i$, denote
$$
V_i := \Phi_i^{-1}(]-1,1[^d\times\left\{t_i\right\})/
$$
By definition,
$$
Y\subseteq\munion_{i=1}^mV_i.
$$
For all $i$, upon reducing $U_i$ slightly if necessary, we may suppose that $V_i$ is relatively compact as a subset of $\Sigma$. By Lemma \procref{SinglePlaque} again, upon reducing each $U_i$ further if necessary, we may suppose that, for all $i$ and for all $j$,
$$
V_i\minter U_j = V_i\minter V_j.
$$
For all $i\in\left\{1,...,m\right\}$, let
$$\eqalign{
p_{i,1}:]-1,1[^d\times T_i&\rightarrow]-1,1[^d\ \text{and}\cr
p_{i,2}:]-1,1[^d\times T_i&\rightarrow T_i\cr}
$$
be the projections onto the first and second factors respectively. For all $i$, define $\pi_i:U_i\rightarrow\Sigma$ by
$$
\pi_i(x) = \Phi_i^{-1}(p_{i,1}(\Phi_i(x)),t_i).
$$
Observe that, for all $i$, $\pi_i$ is leafwise smooth and, for all $y\in V_i$, $\pi_i(y)=y$.
\par
The open set $\Omega$ and the function $\pi$ will be constructed by induction as follows. For all $i$, define
$$
W_i := V_1\munion...\munion V_i.
$$
Suppose that, for some $i$, we have constructed an open set $\Omega_i$ and a leafwise smooth map $\tilde{\pi}_i:\Omega_i\rightarrow\Sigma$ such that $W_i\subseteq\Omega_i$ and $\tilde{\pi}_i(y)=y$ for all $y\in W_i$. Let $O$ be a neighbourhood of the diagonal in $\Sigma\times\Sigma$ consisting of pairs $(p,q)$ for which there exists a unique length-minimising geodesic $\gamma_{pq}:[0,1]\rightarrow\Sigma$ such that $\gamma_{pq}(0)=p$ and $\gamma_{pq}(1)=q$. Define the smooth map $G:O\times[0,1]\rightarrow\Sigma$ by
$$
G(p,q,t) := \gamma_{pq}(t).
$$
Upon reducing $\Omega_i$ and $U_{i+1}$ if necessary, we may suppose that $\Omega_i=\tilde{\pi}_i^{-1}(W_i)$ and that $U_{i+1}=\pi_{i+1}^{-1}(V_{i+1})$. By Theorem \procref{PartitionOfUnity}, there exists a leafwise smooth partition of unity $(f_i,g_i)$ of $\Omega_i\munion U_{i+1}$ subordinate to the cover $(\Omega_i,U_{i+1})$. Define the open subsets $\Omega_{i+1,1}\subseteq\Omega_i$, $\Omega_{i+1,2}\subseteq U_{i+1}$ and $\Omega_{i+1,3}\subseteq\Omega_i\minter U_{i+1}$ by
$$\eqalign{
\Omega_{i+1,1} &:= \Omega_i\setminus\opSupp(g_i),\cr
\Omega_{i+1,2} &:= U_{i+1}\setminus\opSupp(f_i),\ \text{and}\cr
\Omega_{i+1,3} &:= (\tilde{\pi}_i,\pi_{i+1})^{-1}(O).\cr}
$$
Denote
$$
\Omega_{i+1} := \Omega_{i+1,1}\munion\Omega_{i+1,2}\munion\Omega_{i+1,3}.
$$
We claim that $W_{i+1}:=W_i\munion V_{i+1}\subseteq\Omega_{i+1}$. Indeed, $W_i\setminus\opSupp(g)\subseteq\Omega_{i+1,1}$ and $V_{i+1}\setminus\opSupp(f)\subseteq\Omega_{i+1,2}$. By construction,
$$
W_i\minter\opSupp(g)\subseteq W_i\minter U_{i+1}\subseteq V_{i+1},
$$
so that $\tilde{\pi}$ and $\pi_{i+1}$ are both defined and are equal over this set. Likewise,
$$
V_{i+1}\minter\opSupp(f)\subseteq V_{i+1}\minter\big(\munion_{j=1}^iU_j\big) = V_{i+1}\minter\big(\munion_{j=1}^i V_j\big) \subseteq W_i,
$$
so that $\tilde{\pi}$ and $\pi_{i+1}$ again are both defined and are equal over this set. It follows that
$$
W_i\minter\opSupp(g),V_{i+1}\minter\opSupp(f)\subseteq\Omega_{i+1,3},
$$
so that $W_i\munion V_{i+1}\subseteq\Omega_{i+1}$, as asserted.
\par
Define $\pi_{i+1}:\Omega_{i+1}\rightarrow\Sigma$ by
$$
\pi_{i+1}(x):=
\left\{
\matrix\tilde{\pi}_i(x)\hfill&\ \text{if}\ x\in \Omega_{i+1,1},\hfill\cr
\pi_{i+1}(x)\hfill&\ \text{if}\ x\in \Omega_{i+1,2},\ \text{and}\hfill\cr
G(\tilde{\pi}_i(x),\pi_{i+1}(x),1-f_i(x))\hfill&\ \text{if}\ x\in \Omega_{i+1,3}.\hfill\cr\endmatrix\right.
$$
We readily verify that $\pi_{i+1}$ is well-defined and, for all $y\in W_{i+1}$, $\tilde{\pi}_{i+1}(y)=y$. Furthermore, since this function is leafwise smooth over each of $\Omega_{i+1,1}$, $\Omega_{i+1,2}$ and $\Omega_{i+1,3}$, it is leafwise smooth over the whole of $\Omega_{i+1}$. Furthermore, $W_{i+1}\subseteq\Omega_{i+1}$ and,
\par
The induction process is initiated with $W_1:=V_1$, $\Omega_1:=U_1$ and $\tilde{\pi}_1:=\pi_1$ and the result follows upon setting $\Omega:=\Omega_m$ and $\pi:=\tilde{\pi}_m$.\qed
\proclaim{Lemma \nextprocno}
\noindent Let $\Omega\subseteq Y\subseteq X$ be such that $\Omega$ is open in the leafwise topology and $Y$ is compact in the leafwise topology. There exists an open subset $\hat{\Omega}$ of $X$ such that $\hat{\Omega}\minter Y=\Omega$.
\endproclaim
\proclabel{LocalExtensionOfOpenSets}
\proof Let $x$ be a point of $\Omega$. Let $(U,T,\Phi)$ be a laminated chart of $X$ containing $x$. By Lemma \procref{SinglePlaque}, we may suppose that the only plaque of $U$ which intersects $Y$ non-trivially is the plaque containing $x$. Upon reducing $U$ if necessary, we may suppose furthermore that this plaque is contained entirely in $\Omega$. The result follows upon taking the union of all such laminated charts.\qed
\proclaim{Lemma \nextprocno}
\noindent Let $Y\subseteq\Omega\subseteq X$ be such that $Y$ is compact in the leafwise topology and $\Omega$ is open in the leafwise topology. There exists a sequence $(C_k)_{k\in\Bbb{N}}$ of positive constants such that if $f\in C_0^\infty(\Omega)$ is supported in $Y$ and if $\hat{\Omega}$ is an open subset of $X$ containing $\Omega$, then there exists a leafwise smooth function $\tilde{f}\in C^\infty_l(X)$ supported in $\hat{\Omega}$ such that
\medskip
\myitem{(1)} for all $x\in Y$, $\tilde{f}(x)=f(x)$; and
\medskip
\myitem{(2)} for all $k$,
$$
\|\tilde{f}\|_{C^k_l} \leq C_k \|f\|_{C^k}.
$$
\endproclaim
\proclabel{ExtensionOfBumpFunctions}
\proof Suppose first that $Y$ is contained in a single plaque of a laminated chart $(U,T,\Phi)$ and that this plaque is contained in $\Omega$. By Lemma \procref{SinglePlaque}, we may suppose that
$$
\Phi(Y)\subseteq ]-1,1[^d\times\left\{t\right\},
$$
for some $t\in T$. Let
$$\eqalign{
p_1:]-1,1[^d\times T&\rightarrow]-1,1[^d\ \text{and}\cr
p_2:]-1,1[^d\times T&\rightarrow T\cr}
$$
be the projections onto the first and second factors respecively. Define $\pi:U\rightarrow\Sigma$ by
$$
\pi(x) := \Phi^{-1}(p_1(\Phi(x)),t).
$$
Let $d$ be a metric of $T$. Let $\chi\in C^\infty(\Bbb{R},[0,1])$ be such that $\chi(t)=1$ for $t\leq 1$ and $\chi(t)=0$ for $t\geq 2$. For $\lambda>0$, define $\chi_\lambda:U\rightarrow\Bbb{R}$ by
$$
\chi_\lambda(x) := \chi(d(p_2(\Phi(x)),t)/\lambda).
$$
Define $\tilde{f}_\lambda$ by
$$
\tilde{f}_\lambda(x) := (f\circ\pi)(x)\chi_\lambda(x).
$$
\par
Since $T$ is locally compact, there exists $\epsilon>0$ such that $Z:=\Phi^{-1}(\Phi(Y)\times\overline{B}_\epsilon(t))$ is a compact, subset of $X$, where $B_\epsilon(t)$ denotes the ball of radius $\epsilon$ about $t$ in $T$ with respect to the metric $d$. Upon reducing $\epsilon$ if necessary, we may suppose in addition that $Z\subseteq\hat{\Omega}$. For sufficiently large $\lambda$, $\opSupp(\tilde{f}_\lambda)\subseteq Z\subseteq\hat{\Omega}$. Furthermore, by compactness, for all $k$, there exists $C_k$ such that, for all such $\lambda$,
$$
\|\tilde{f}_\lambda\|_{C^k_l} \leq \|f\circ\pi|_Z\|_{C^k_l} \leq C_k \|f\|_{C^k}.
$$
Finally, by compactness again, $Y$ is contained in finitely many such laminated charts, and the general case follows using a finite $C^\infty_l$ partition of unity.\qed
\newsubhead{Hyperbolic perturbations}[HyperbolicPerturbations]
Before proving our main result concerning the infinite dimensionality of Teichm\"uller space, we return to the case of a single complete, hyperbolic surface $\Sigma$ without cusps. Let $g$ denote its hyperbolic metric. We review the perturbation theory of hyperbolic metrics over this surface, adopting a formalism similar to that of \cite{Tromba}. Let $\Sigma$ be a complete hyperbolic surface without cusps. Let $A$ be a smooth, bounded section of $\opEnd(T\Sigma)$ which is symmetric with respect to $g$. For sufficiently small $t$, denote
$$
g_t := g((\opId + tA)\cdot,(\opId + tA)\cdot),\eqnum{\nexteqnno[VariationOfMetricsFromA]}
$$
and let $\kappa_t:\Sigma\rightarrow\Bbb{R}$ be the curvature function of this metric. $A$ is said to be a {\sl hyperbolic perturbation} whenever
$$
\kappa_t = -1 + \opO(t^2).\eqnum{\nexteqnno[HyperbolicPerturbationDefinition]}
$$
Given a hyperbolic perturbation $A$, for every simple, closed geodesic $\gamma$ in $\Sigma$, define
$$
\Delta([\gamma],A) := \left.\frac{\partial}{\partial t}\opLog(l([\gamma],h_t))\right|_{t=0},
$$
where $l([\gamma],h_t)$ is as in Section \subheadref{MainResults}. We readily obtain
\proclaim{Lemma \nextprocno}
\noindent If $\gamma$ is parametrised by arc length, then
$$
\Delta([\gamma],A) = \frac{1}{l([\gamma],g)}\int_\gamma\langle\dot{\gamma},A\dot{\gamma}\rangle dl.
$$
\endproclaim
\proclabel{FormulaForDelta}
Since the curvature operator is a second order, non-linear partial differential operator, \eqnref{HyperbolicPerturbationDefinition} can be rewritten as a linear differential condition on $A$. In particular, the space of hyperbolic perturbations is a vector space. It will be helpful for what follows to determine this condition explicitly. To this end, recall first that the {\sl divergence} of an endomorphism field $B$ with respect to $g$ is the $1$-form defined by
$$
(\nabla\cdot B)(\cdot) = \sum_{i=1}^2\langle e_i,(\nabla_{e_i}B)(\cdot)\rangle,\eqnum{\nexteqnno[FirstDefinitionOfDivergence]}
$$
where $(e_1,e_2)$ is a local orthonormal frame of $g$. Likewise, the {\sl divergence} of a $1$-form $\beta$ with respect to $g$ is the function defined by
$$
\nabla\cdot\beta = \sum_{i=1}^2(\nabla_{e_i}\beta)(e_i).\eqnum{\nexteqnno[SecondDefintionOfDivergence]}
$$
\proclaim{Lemma \nextprocno}
\noindent Let $A$ be a smooth, bounded, symmetric section of $\opEnd(T\Sigma)$. $A$ is a hyperbolic perturbation if and only if
$$
\nabla\cdot\nabla\cdot(JAJ) + \opTr(A) = 0,\eqnum{\nexteqnno[HyperbolicPerturbationCondition]}
$$
where $J$ denotes the complex structure of $g$.
\endproclaim
\proclabel{HyperbolicPerturbationCondition}
\proof As in the proof of Lemma \procref{PerturbedCurvatureFormula}, we use the formalism of moving frames. Let $(e_1,e_2)$ denote a local orthonormal moving frame of $g$. For all sufficiently small $t$, denote $B_t:=\opId+tA$ and
$$
g_t := g(B_t\cdot,B_t\cdot),
$$
and let $\kappa_t$ denote the curvature of $g_t$. For all such $t$, a locally orthonormal frame of $g_t$ is given by $(B_t^{-1}e_1,B_t^{-1}e_2)$. Using \eqnref{ConnectionFormFormula}, we show that the connection form of this frame is
$$
\alpha_t(X) = \alpha_0(X) - \opDet(B_t^{-1})g(B_t\nabla\cdot(B_tJ),X).
$$
Differentiating at $t=0$ then yields
$$
\left.\frac{\partial}{\partial t}\alpha_t(X)\right|_{t=0} = -g(\nabla\cdot(AJ),X),
$$
so that
$$
\left.\frac{\partial}{\partial t}d\alpha_t(e_1,e_2)\right|_{t=0} = \nabla\cdot\nabla\cdot(JAJ).
$$
It follows by \eqnref{MovingFrameCurvature} that
$$
\left.\frac{\partial}{\partial t}\kappa_t\right|_{t=0} = \nabla\cdot\nabla\cdot(JAJ) - \kappa_g\opTr(A),
$$
and the result follows since $\kappa_g=-1$.\qed
\newsubhead{Infinite dimensionality of Teichm\"uller space}[InfiniteDimensionalityOfTeichmuellerSpace]
Let $X$ be a compact hyperbolic surface lamination and denote its leafwise metric by $g$. By compactness, the injectivity radius of every leaf is bounded below by $\epsilon>0$, say, so that no leaf of $X$ has cusps. Recall (c.f. \cite{SullivanI} and \cite{SullivanII}) that the Teichm\"uller space of $X$ is defined to be the space of transversally continuous conformal structures over $X$ modulo leafwise diffeomorphisms of $X$ which are leafwise isotopic to the identity. Recall, furthermore, that \cite{SullivanII}, Sullivan shows that this space naturally carries the structure of a Banach manifold. In this section, we show
\proclaim{Theorem \nextprocno}
\noindent Suppose that, for every $K\subseteq X$ which is compact in the leafwise topology, there exists a simple, closed geodesic $\gamma$ in $X$ such that
$$
\gamma\minter K = \emptyset.
$$
Then there exists a sequence $(\gamma_m)$ of simple, closed geodesics in $X$ such that, for every finite sequence $a_1,...,a_m\in\Bbb{R}$, there exists a smooth family $(g_t)_{t\in]-\epsilon,\epsilon[}$ of leafwise smooth hyperbolic metrics such that, for all $1\leq i\leq m$,
$$
\left.\frac{\partial}{\partial t}\opLog(l([\gamma_i],g_t))\right|_{t=0} = a_i.
$$
In particular, the Teichm\"uller space of $X$ is infinite dimensional.
\endproclaim
\proclabel{InfiniteDimensionality}
Theorem \procref{InfiniteDimensionality} is proven by a straightforward induction argument. The induction step is provided by
\proclaim{Lemma \nextprocno}
\noindent Suppose that for every $K\subseteq X$ which is compact in the leafwise topology, there exists a simple, closed geodesic $\gamma$ in $X$ such that
$$
\gamma\minter K = \emptyset.
$$
For every finite set $G$ of simple, closed geodesics in $X$, and for all $\epsilon>0$, there exists a simple, closed geodesic $\gamma$ in $X$ and a hyperbolic perturbation $A$ such that
$$
\Delta([\gamma],A) < -\frac{1}{2\pi},
$$
and, for all $\gamma'\in G$,
$$
\left|\Delta([\gamma'],A)\right| < \epsilon.
$$
\endproclaim
\proclabel{InductionStep}
\proof Let $R_1>0$ be a positive number to be determined presently. Let $\gamma$ be a simple, closed geodesic in $X$ such that, for all $\gamma'\in G$,
$$
d_l(\gamma,\gamma') > 2R_1.
$$
Let $\Sigma$ be the leaf containing $\gamma$. Let $G_0$ be the subset of $G$ consisting of those geodesics that are contained in $\Sigma$ and let $G_1$ consist of all other geodesics of $G$.
\par
Let $R<R_1$ be such that the tubular neighbourhood $B_R(\gamma)$ of radius $R$ about $\gamma$ is isometric to $\gamma\times]-R,R[$ furnished with the twisted product metric
$$
\opCosh^2(r)dt^2 + dr^2.
$$
Let $\psi\in C^\infty_0(]-R/2,R/2[)$ be a positive function such that
$$
\int_{-R}^R\frac{\psi(r)}{\opCosh(r)}dr = 1.
$$
Define the endomorphism field $A_1$ over $B_R(\gamma)$ by
$$
g(A_1\cdot,\cdot) = \psi(r)dr^2,
$$
and extend it to a smooth endomorphism field over the whole of $\Sigma$ by setting it equal to zero outside this annulus. Let $u_1:\Sigma\rightarrow\Bbb{R}$ be the unique bounded, smooth function such that
$$
(\Delta - 2)u_1 = \nabla\cdot\nabla\cdot(JA_1J) + \opTr(A_1),
$$
so that
$$
A_2 := A_1 + u_1\opId,
$$
is a hyperbolic perturbation. By Theorem \procref{RatesOfChangeOfLength} and Lemma \procref{ModulesOfSmoothTransformations},
$$
\Delta([\gamma],A_2) \leq -\frac{1}{\pi},
$$
and, provided $R_1>0$ is sufficiently large, for all $\gamma'\in G_0$,
$$
\left|\Delta([\gamma'],A_2)\right|<\frac{\epsilon}{2}.
$$
\par
Let $R_2>R_1$ be another positive constant to be determined presently. Let $\chi\in C_0^\infty(B_{3R/4}(\gamma),[0,1])$ be such that $\chi$ is equal to $1$ over $\overline{B}_{R/2}(\gamma)$. By Lemma \procref{TubularNeighbourhood}, there exists a neighbourhood $\Omega_1$ of $\overline{B}_R(\gamma)$ in $X$ and a leafwise smooth function $\pi:\Omega_1\rightarrow\Sigma$ such that, for all $x\in\overline{B}_R(\gamma)$, $\pi(x)=x$. By Lemma \procref{LocalExtensionOfOpenSets}, there exists another open subset $\Omega_2$ of $X$ such that
$$
\Omega_2\minter\big(B_{R_2}(\gamma)\munion\munion_{\gamma'\in G}B_{R_2}(\gamma')\big) = B_R(\gamma).
$$
Set $\Omega:=\Omega_1\minter\Omega_2$. By Lemma \procref{ExtensionOfBumpFunctions}, there exists $C_1>0$, independent of $R_2$ and $\Omega$, and a leafwise smooth function $\tilde{\chi}:X\rightarrow[0,1]$ such that
\medskip
\myitem{(1)} $\opSupp(\tilde{\chi})\subseteq\Omega$,
\medskip
\myitem{(2)} $\tilde{\chi}(x)=\chi(x)$ for all $x\in B_R(\gamma)$, and
\medskip
\myitem{(3)} $\|\tilde{\chi}\|_{C^2_l}\leq C_1$.
\medskip
Define the leafwise smooth endomorphism field $A_3$ over $X$ by
$$
A_3 :=\left\{\matrix\tilde{\chi}(x)(\pi^* A_1)(x)\hfill&\ \text{if}\ x\in\Omega\ \text{and}\hfill\cr
0\hfill&\ \text{if}\ x\notin\opSupp(\tilde{\chi})\hfill\cr\endmatrix\right.,
$$
and let $u_3:\Sigma\rightarrow\Bbb{R}$ be the unique bounded, leafwise smooth function such that
$$
Lu_3 = \nabla\cdot\nabla\cdot(JA_3J) + \opTr(A_3),
$$
so that
$$
A_4 := A_3 + u_3\opId,
$$
is a hyperbolic perturbation. For $\gamma'\in G_0\munion\left\{\gamma\right\}$, we have
$$\eqalign{
\opSupp(A_3-A_1)\minter B_{R_2}(\gamma') &= \emptyset\ \text{and}\cr
\opSupp\big((\Delta-2)(u_3-u_1)\big)\minter B_{R_2}(\gamma') &= \emptyset.\cr}
$$
It follows from the first relation that, near $\gamma'$, $A_4-A_2$ only depends on $u_3-u_1$. It then follows by Lemmas \procref{ExponentialDecayAwayFromSource} and \procref{FormulaForDelta} that there exists $C_2>0$, independent of $R_2$, such that
$$
\left|\Delta([\gamma'],A_4) - \Delta([\gamma'],A_2)\right| \leq C_2 e^{-R_2}.
$$
Likewise, for $\gamma'\in G_1$,
$$\eqalign{
\opSupp(A_3)\minter B_{R_2}(\gamma') &= \emptyset\ \text{and}\cr
\opSupp\big(L(u_3)\big)\minter B_{R_2}(\gamma') &= \emptyset,\cr}
$$
so that, upon increasing $C_2$ if necessary, we have
$$
\left|\Delta([\gamma'],A_4)\right| \leq C_2e^{-R_2}.
$$
The result now follows upon setting $R_2$ sufficiently large.\qed
\proclaim{Lemma \nextprocno}
\noindent For every leafwise smooth hyperbolic perturbation $A$ over $X$ there exists $\epsilon>0$ and a smooth family $(g_t)_{t\in]-\epsilon,\epsilon[}$ of leafwise metrics over $X$ such that
$$
\left.\frac{\partial}{\partial t}g_t\right|_{t=0} = 2g(A\cdot,\cdot).
$$
\endproclaim
\proclabel{IntegrabilityOfHyperbolicPerturbations}
\proof Indeed, for sufficiently small $t$, denote $B_t:=\opId + tA$ and
$$
h_t := g(B_t\cdot,B_t\cdot).
$$
By Lemma \procref{CandelSmoothness}, there exists $\epsilon>0$ and a smooth family $(u_t)_{t\in]-\epsilon,\epsilon[}$ of leafwise smooth functions such that, for all $t$, the metric
$$
g_t := e^{2u_t}h_t
$$
is hyperbolic. Consider now a leaf $\Sigma$ of $X$. By Lemma \procref{HyperbolicPerturbationCondition}, over $\Sigma$,
$$
L\left.\frac{\partial}{\partial t}u_t\right|_{t=0} = \nabla\cdot\nabla\cdot(JAJ) + \opTr(A) = 0.
$$
However, since $\Sigma$ is without cusps, $L$ has trivial kernel in the space $C^2_\opbdd(\Sigma)$ of bounded, twice differentiable functions over this surface, so that, over $\Sigma$,
$$
\left.\frac{\partial}{\partial t}u_t\right|_{t=0} = 0.
$$
Since $\Sigma$ is arbitrary, this relation in fact holds over the whole of $X$, so that, since $u_0=0$,
$$
\left.\frac{\partial}{\partial t}g_t\right|_{t=0} = \left.\frac{\partial}{\partial t}h_t\right|_{t=0} = 2g(A\cdot,\cdot),
$$
as desired.\qed
\medskip
{\bf\noindent Proof of Theorem \procref{InfiniteDimensionality}:\ }Using Lemma \procref{InductionStep}, we show by induction that there exist sequences $(\epsilon_m)$, $(B_m)$ of positive numbers, a sequence $(\gamma_m)$ of simple, closed geodesics in $X$, and a sequence $(A_m)$ of hyperbolic perturbations such that
\medskip
\myitem{(1)} for all $m$,
$$
\epsilon_{m+1}\leq \frac{1}{(4\pi)^{m+1}(m+1)!B_1\cdots B_m};
$$
\myitem{(2)} for all $m$,
$$
\Delta([\gamma_m],A_m)\leq -\frac{1}{2\pi};
$$
\myitem{(3)} for all $m$, and for all $n<m$,
$$
\left|\Delta([\gamma_n],A_m)\right| < \epsilon_m;\ \text{and}
$$
\myitem{(4)} for all $m$, and for every simple, closed geodesic $\gamma$,
$$
\left|\Delta([\gamma],A_m)\right| < B_m.
$$
Consider now a finite sequence $a_1,...,a_m\in\Bbb{R}$, and define the $m\times m$ matrix $M$ by
$$
M_{ij} := \Delta([\gamma_i],A_j).
$$
The above estimates yield,
$$
\left|\opDet(M)\right| \geq \frac{1}{(4\pi)^m},
$$
so that $M$ is invertible. There therefore exists $\alpha_1,...,\alpha_m\in\Bbb{R}$ such that, for all $1\leq i\leq m$,
$$
\Delta\big([\gamma_i],\sum_{j=1}^m\alpha_j A_j\big) = a_i.
$$
Setting
$$
A:=\sum_{j=1}^m\alpha_j A_j,
$$
the result follows by Lemma \procref{IntegrabilityOfHyperbolicPerturbations}.\qed
\newhead{Bibliography}[Bibliography]
\medskip
{\leftskip = 5ex \parindent = -5ex
\leavevmode\hbox to 4ex{\hfil \cite{Ahlfors}}\hskip 1ex{Ahlfors L. V., {\sl Conformal Invariants: Topics in Geometric Function Theory}, McGraw-Hill, New York, D\"usseldorf, Johannesburg, (1973)}
\medskip\leavevmode\hbox to 4ex{\hfil \cite{AlcaldeCuestaDalBoMartinezVerjovsky}}\hskip 1ex{Alcalde Cuesta F., Dal'Bo F., Mart\'\i nez M., Verjovsky, Minimality of the horocycle flow
on laminations by hyperbolic surfaces with non-trivial topology, {\sl Discrete Contin. Dyn. Syst.}, {\bf 36}, no. 9, (2016), 4619--4635}
\medskip\leavevmode\hbox to 4ex{\hfil \cite{AlessandriniLiuPapadopoulosSuI}}\hskip 1ex{Alessandrini A., Liu L., Papadopoulos A., Su W., Sun Z., On Fenchel-Nielsen coordinates
on Teichm\"ller spaces of surfaces of infinite type, {\sl Ann. Acad. Sci. Fenn. Math.}, {\bf 36}, no. 2, (2011), 621--659}
\medskip\leavevmode\hbox to 4ex{\hfil \cite{AlessandriniLiuPapadopoulosSuII}}\hskip 1ex{Alessandrini A., Liu L., Papadopoulos A., Su W., On local comparison between various
metrics on Teichm\"uller spaces, {\sl Geom. Dedicata}, {\bf 157}, (2012), 91--110}
\medskip\leavevmode\hbox to 4ex{\hfil \cite{AlessandriniLiuPapadopoulosSuIII}}\hskip 1ex{Alessandrini A., Liu L., Papadopoulos A., Su W., On various Teichm\"uller spaces of a
surface of infinite topological type, {\sl Proc. Amer. Math. Soc.}, {\bf 140}, no. 2, (2012), 561--574}
\medskip\leavevmode\hbox to 4ex{\hfil \cite{AlvarezBrum}}\hskip 1ex{Alvarez S., Brum J., Topology of leaves for minimal laminations by hyperbolic surfaces II, Preprint, 2019}
\medskip\leavevmode\hbox to 4ex{\hfil \cite{AlvarezBrumMartinezPotrie}}\hskip 1ex{Alvarez S., Brum J., Mart\'\i nez, Potrie R., Topology of leaves for minimal laminations by hyperbolic surfaces, with an appendix written with M. Wolff. Preprint, 2019}
\medskip\leavevmode\hbox to 4ex{\hfil \cite{AlvarezLessa}}\hskip 1ex{Alvarez S., Lessa P., The Teichmüller space of the Hirsch foliation, {\sl Ann. Inst. Fourier}, {\bf 68}, no. 1, (2018), 1--51}
\medskip\leavevmode\hbox to 4ex{\hfil \cite{BallmanSchroederGromov}}\hskip 1ex{Ballmann W., Gromov M., Schroeder V., {\sl Manifolds of Nonpositive Curvature},\break Progress in Mathematics, {\bf 61}, Birkha\"user Verlag, (1985)}
\medskip\leavevmode\hbox to 4ex{\hfil \cite{CandelI}}\hskip 1ex{Candel A., Uniformization of surface laminations, {\sl Ann. Sci. ENS.}, {\bf 4}, no. 26, (1993), 489--516}
\medskip\leavevmode\hbox to 4ex{\hfil \cite{CandelConlon}}\hskip 1ex{Candel A., Conlon L., {\sl Foliations I}, Graduate Studies in Mathematics, {\bf 23}, AMS, Providence, Rhode Island, (2000)}
\medskip\leavevmode\hbox to 4ex{\hfil \cite{Chavel}}\hskip 1ex{Chavel I., {\sl Riemannian Geometry: A Modern Introduction}, Cambridge Studies in Advanced Mathematics, {\bf 98}, CUP, (2006)}
\medskip\leavevmode\hbox to 4ex{\hfil \cite{ChowKnopf}}\hskip 1ex{Chow B., Knopf D., {\sl The Ricci Flow: An Introduction}, Mathematical Surveys and Monographs, {\bf 110}, AMS, (2004)}
\medskip\leavevmode\hbox to 4ex{\hfil \cite{Deroin}}\hskip 1ex{Deroin B., Nonrigidity of hyperbolic surfaces laminations, {\sl Proc. Amer. Math. Soc.}, {\bf 135}, no. 3, (2007), 873--881}
\medskip\leavevmode\hbox to 4ex{\hfil \cite{DumasWolf}}\hskip 1ex{Dumas D., Wolf M., Projective structures, grafting and measured laminations, {\sl Geom. Topol.}, {\bf 12}, no. 1, (2008), 351--386}
\medskip\leavevmode\hbox to 4ex{\hfil \cite{EpsteinMillettTischler}}\hskip 1ex{Epstein D. B. A., Millett K. C., Tischler D., Leaves without holonomy, {\sl J. London Math. Soc.}, {\bf 16}, no. 3, (1977), 548--552}
\medskip\leavevmode\hbox to 4ex{\hfil \cite{GhysI}}\hskip 1ex{Ghys \'E., Laminations par surfaces de Riemann, in {\sl Dynamique et g\'eom\'trie complexes (Lyon,
1997)}, {\sl Panor. Synth\`eses}, Volume 8, Pages ix, xi, 49--95, Soc. Math. France, Paris, 1999}
\medskip\leavevmode\hbox to 4ex{\hfil \cite{GilbTrud}}\hskip 1ex{Gilbarg D., Trudinger N., {\sl Elliptic partial differential equations of second order},\break Grundlehren der Mathematischen Wissenschaften, {\bf 224}, Springer-Verlag, Berlin, second edition, 1983}
\medskip\leavevmode\hbox to 4ex{\hfil \cite{Hector}}\hskip 1ex{Hector G., Feuilletages en cylindres, in {\sl Geometry and topology (Proc. III Latin Amer. School
of Math., Inst. Mat. Pura Aplicada CNPq, Rio de Janeiro, 1976)}, Volume 597 of {\sl Lecture Notes in Math.}, Pages 252--270, Springer, Berlin, (1977)}
\medskip\leavevmode\hbox to 4ex{\hfil \cite{HirschPughShub}}\hskip 1ex{Hirsch M. W., Pugh C. C., Shub M., {\sl Invariant manifolds}, Lecture Notes in Mathematics, {\bf 583}, Springer Verlag, Berling, New York, (1977)}
\medskip\leavevmode\hbox to 4ex{\hfil \cite{KatokHasselblatt}}\hskip 1ex{Katok A., Hasselblatt B., {\sl Introduction to the modern theory of dynamical systems}, Encylopedia of Mathematics and its Applications, {\bf 54}, CUP, Cambridge, (1995)}
\medskip\leavevmode\hbox to 4ex{\hfil \cite{LyubichMinsky}}\hskip 1ex{Lyubich M., Minsky Y., Laminations in holomorphic dynamics, {\sl J. Diff. Geom.}, {\bf 47}, (1997), 17--94}
\medskip\leavevmode\hbox to 4ex{\hfil \cite{McMullen}}\hskip 1ex{McMullen C., Complex earthquakes and Teichmüller theory, {\sl J. Amer. Math. Soc.}, {\bf 11}, no. 2, (1998), 283--320}
\medskip\leavevmode\hbox to 4ex{\hfil \cite{MooreSchochet}}\hskip 1ex{Moore C. C., Schochet C. L., {\sl Global analysis on foliated spaces}, MSRI Publications, {\bf 9}, Cambridge University Press, New York, second edition, 2006}
\medskip\leavevmode\hbox to 4ex{\hfil \cite{Rosenberg}}\hskip 1ex{Rosenberg H., Foliations by planes, {\sl Topology}, {\bf 7}, (1968), 131--138}
\medskip\leavevmode\hbox to 4ex{\hfil \cite{Saric}}\hskip 1ex{\v Sari\'c, The Teichm\"uller theory of the solenoid, in {\sl Handbook of Teichmüller theory. Vol. II},
Volume 13 of {\sl IRMA Lect. Math. Theor. Phys.}, Pages 811--857, Eur. Math. Soc., Zürich, 2009}
\medskip\leavevmode\hbox to 4ex{\hfil \cite{ScannellWolf}}\hskip 1ex{Scannell K., Wolf M., The grafting map of Teichmüller space, {\sl J. Amer. Math. Soc.}, {\bf 15}, no. 4, (2002), 893--927}
\medskip\leavevmode\hbox to 4ex{\hfil \cite{Spivak}}\hskip 1ex{Spivak M., {\sl A Comprehensive Introduction to Differential Geometry}, Vols. 1-5, Publish of Perish, (1999)}
\medskip\leavevmode\hbox to 4ex{\hfil \cite{SullivanI}}\hskip 1ex{Sullivan D., Bounds, quadratic differentials, and renormalization conjectures, in\break {\sl American Mathematical Society centennial publications, Vol. II, (Providence, RI, 1988)}, Pages 417--466, Amer. Math. Soc., Providence, RI, 1992}
\medskip\leavevmode\hbox to 4ex{\hfil \cite{SullivanII}}\hskip 1ex{Sullivan D, Linking the universalities of Milnor-Thurston, Feigenbaum and Ahlfors-Bers, in {\sl Topological methods in modern mathematics (Stony Brook, NY, 1991)}, Pages 543--564, Publish or Perish, Houston, TX, 1993}
\medskip\leavevmode\hbox to 4ex{\hfil \cite{Tromba}}\hskip 1ex{Tromba A., {\sl Teichm\"uller theory in Riemannian Geometry}, Lectures in Mathematics, ETH Z\"urich, Birkha\"user Verlag, 1992}
\medskip\leavevmode\hbox to 4ex{\hfil \cite{Verjovsky}}\hskip 1ex{Verjovsky A., A uniformization theorem for holomorphic foliations, in {\sl The Lefschetz centennial conference, Part III (Mexico City, 1984)}, Volume 58 of {\sl Contemp. Math.}, Pages 233--253, Amer. Math. Soc., Providence, RI, 1987}
\medskip\leavevmode\hbox to 4ex{\hfil \cite{Wolpert}}\hskip 1ex{Wolpert S., An elementary formula for the Fenchel-Nielsen twist, {\sl Comment. Math. Helv.}, {\bf 56}, no. 1, (1981), 132--135}
\par}
%
%%%%%%%%%%%%%%%%%%%%%%%%%%%%%%%%%%%%%%%%%%%%%%%%%%%%%%%%%%%%%%%%%%%%%%%%%%%%%%%%%%%%%%%%%%%%%%%%%%%%%%%%%%%%%%%%%%%%%%%
%
% 3: Closing commands.
%
%%%%%%%%%%%%%%%%%%%%%%%%%%%%%%%%%%%%%%%%%%%%%%%%%%%%%%%%%%%%%%%%%%%%%%%%%%%%%%%%%%%%%%%%%%%%%%%%%%%%%%%%%%%%%%%%%%%%%%%
%
\enddocument